\documentclass{amsart}
\usepackage{amssymb}
\newcommand\datver[1]{\def\datverp
 {\par\boxed{\boxed{\text{Version: #1; Run: \today}}}}}
\datver{3.0B; Revised: 6/11/2001}
\newcommand\CC{\mathbb C}

\newcommand\RR{\mathbb R}
\newcommand\ZZ{\mathbb Z}

\newcommand\pa{\partial}
\newcommand\GR{\mathcal{G}}
\newcommand\CI{\mathcal{C}^\infty}
\newcommand\CIc{\mathcal{C}_c^\infty}
\newcommand\Tt{a_0 \otimes a_1 \otimes \ldots \otimes a_n}
\newcommand\ie{i.e. }
\newcommand\Diff{\operatorname{Diff}}
\newcommand\Diffeo{\operatorname{Diffeo}}
\newcommand\B{\operatorname{B}}
\newcommand\II{\operatorname{I}}
\newcommand\End{\operatorname{End}}
\newcommand\Exp{\operatorname{Exp}}
\newcommand\HH{\operatorname{HH}}
\newcommand\Hd{\operatorname{HH}}
\newcommand\Hc{\operatorname{HC}}
\newcommand\HC{\operatorname{HC}}
\newcommand\HP{\operatorname{HP}}
\newcommand\Hp{\operatorname{HP}}
\newcommand\EH{\operatorname{EH}}
\newcommand\EC{\operatorname{EC}}

\newcommand{\tPS}[1]{\Psi^{#1}(\GR)}
\newcommand{\PS}[1]{\Psi^{#1}(\GR;E)}
\newcommand\mO{{\mathcal O}(M)}

\newcommand\Tr{\operatorname{Tr}}
\newcommand\STr{\operatorname{STr}}
\newcommand{\R}{\RR}

\newcommand{\Z}{\ZZ}

\newcommand{\cohom}{\operatorname{H}}
\newcommand\alge{\mathcal A}
\newcommand\ideal{\mathcal I}

\newcommand\maL{\mathcal L}
\newcommand\maH{\mathcal H}

\newcommand\maF{\mathcal F}

\newcommand\maP{\mathcal P}

\newcommand\potimes{\widehat{\otimes}}
\newtheorem{theorem}{Theorem}
\newtheorem{proposition}{Proposition}
\newtheorem{corollary}{Corollary}
\newtheorem{lemma}{Lemma}

\theoremstyle{definition}
\newtheorem{definition}{Definition}

\theoremstyle{remark}
\newtheorem{remark}{Remark}

\begin{document}

%%%%%%%%%%%%%%%%%%%%%%%%%%%%%%%%%%%%%%%%%%%%%%%%%%%%%%%
%                                                     %
%   THE MANUSCRIPT BEGINS HERE                        %
%                                                     %
%%%%%%%%%%%%%%%%%%%%%%%%%%%%%%%%%%%%%%%%%%%%%%%%%%%%%%%

%%% TITLE

\title[Homology of families]
{Homology of algebras of families of pseudodifferential operators}

%%% AND THE AUTHORS ARE:

\author[M. Benameur]{Moulay-Tahar Benameur}
\address{Inst. Desargues, Lyon, France}
\email{benameur@desargues.univ-lyon1.fr}
\author[V. Nistor]{Victor Nistor}
\address{Inst. Desargues and Pennsylvania State University,
University Park, PA 16802}
\email{nistor@math.psu.edu}

\thanks{Partially supported by the NSF Young Investigator Award
DMS-9457859, NSF Grant DMS-9971951 and ``collaborative CNRS/NSF
research grant'' DMS-9981251. This manuscript is available from {\bf
http:{\scriptsize//}www.math.psu.edu{\scriptsize/}nistor{\scriptsize/}}.}

\dedicatory\datverp
%\date\datverp

\begin{abstract}
We compute the Hochschild, cyclic, and periodic cyclic homology groups
of algebras of families of Laurent complete symbols on manifolds with
corners.  We show in particular that the spectral sequence associated
with Hochschild homology degenerates at $E^2$ and converges to
Hochschild homology. As a byproduct, we deduce an identification of
the space of residue traces on fibrations by manifolds with
corners. In the process, we prove several general results about
algebras of complete symbols on manifolds with corners.
\end{abstract}
\maketitle \tableofcontents

\section*{Introduction}

Some of the main tools in the applications of Non-commutative geometry
to index theory and other areas of mathematics are the Hochschild and
periodic cyclic homology groups.  Hochschild homology, for example,
can be used to understand the residue trace introduced by Guillemin
and Wodzicki \cite{Guillemin,Wodzicki}. Other higher residue cocycles
appear when studying more complicated singular spaces.

In this paper, we study the Hochschild homology of certain algebras of
complete symbols. Recall that an algebra of complete symbols is the
quotient of the algebra of all pseudodifferential operators by the
ideal of regularizing (or order $-\infty$) operators.  Previously,
results in this direction were obtained in \cite{BenameurNistor,
BrylinskiGetzler, LauterMoroianu, MelroseNistor, Moroianu, Schrohe},
and by Wodzicki \cite{WodzickiU} (unpublished). See also
\cite{SchroheN}.

Our algebras of complete symbols can be obtained as algebras of
complete symbols on differentiable groupoids
\cite{LauterNistor,Monthubert1,NWX}. For this large class of examples,
it has been shown in \cite{BenameurNistor} that the periodic cyclic
homology can be computed, without any further assumption on the
groupoid under consideration, in terms of the Laurent cohomology
spaces of the cosphere bundle of the associated Lie algebroid.  The
Hochschild homology groups of algebras of complete symbols
on differentiable groupoids, however, cannot be described in general in
a simple, uniform way for all differentiable groupoids.  Finding the
right language in which to express these Hochschild homology groups
seems to be a problem in itself--clearly an interesting one.

In the present paper, we are concerned with a class of algebras of
complete symbols that consist of families of pseudodifferential
operators on a fibration by manifolds with corners and with Laurent
singularities at the hyperfaces (see Section~\ref{Sec.OA} for precise
definitions).  Let us consider then a fibration $\pi:M\to B$ of a
manifold with corners over a smooth manifold $B$. Let
$\alge_{\maL}(M\vert B;E)$ be the algebra of complete symbols
associated to the algebra of smooth families of pseudodifferential
operators on the fibers of $\pi$ and with Laurent singularities at the
hyperfaces. (The precise construction of this algebra is done using
groupoids, see Section~\ref{Sec.OA}.) To our surprise, the resulting
algebra turns out not to depend on the groupoid $\GR$ used to define
it, but only on $\pi : M \to B$, as long as the groupoid $\GR$
satisfies the assumptions \eqref{eq.A1} and \eqref{eq.A2} of Section
\ref{Sec.OA}.

Let $S^*_{vert}(M) := (T^*_{vert}M \smallsetminus 0) / \RR_+^*$ denote
the sphere bundle of the vertical cotangent bundle to the fibration
$\pi : M \to B$. The periodic cyclic homology of the algebra
$\alge_{\maL}(M\vert B;E)$ of Laurent vertical complete symbols with
coefficients in a $\Z/2\Z$-graded vector bundle $E$, is given by
Theorem \ref{periodic}:
\begin{equation}
	\HP_j(\alge_{\maL}(M\vert B;E)) \simeq \bigoplus_{k\in \Z}
	\cohom_{c,\maL}^{j+2k}(S^*_{vert}(M) \times S^1), \quad j=0,1,
\end{equation}
where $\cohom_{c,\maL}$ denotes compactly supported Laurent cohomology
(it is the compactly supported cohomology of a certain space
functorially associated to $S^*_{vert}(M)$).

We also compute the Hochschild homology groups of the algebras
$\alge_{\maL}(M\vert B;E)$. For $B$ reduced to a point, this was done
in \cite{BenameurNistor}. For families, the computation requires new
ingredients, and, in particular, we may obtain infinite dimensional
groups. The result is then:
\begin{equation}
	\HH_m(\alge_{\maL}(M\vert B;E)) \simeq \bigoplus_{k+h=m}
	\Omega_c^h(B, \maF^{2p-k}),
\end{equation}
where $\maF^*$ is the local coefficient system over $B$ given by the
Laurent cohomology of the fibers of $S^*_{vert}(M)\times S^1$ and $p$
is the dimension of the fibers of $\pi:M\to B$ (Theorem
\ref{theorem.fib.cor}).  This result leads, in particular, to an
explicit description of the space of residue traces on the algebras of
families of pseudodifferential operators that we consider. For
example, when $M$ is smooth (no corners) and $\pi$ has connected
fibers, we obtain that
\begin{equation}
	\HH^0(\alge(M\vert B;E)) \simeq {\mathcal
	C}^{-\infty}(B) := \CIc(B)',
\end{equation}
that is, that the space of (super)traces on $\alge(M\vert B;E)$
identifies with the space of distributions on $B$.  (Note that in this
case $\alge_{\maL}(M\vert B;E) = \alge(M\vert B;E)$.) The space of
traces in the general case of fibrations $\pi : M \to B$ when $M$ has
corners is also explicitly described in terms of the minimal faces of
the fibers of $\pi$ (Theorem~\ref{theorem.traces}).

Let us now briefly describe the contents of each section. In
Section~\ref{Sec.filtered}, we briefly recall the basic definitions of
$\Z/2\Z$-graded homologies for topologically filtered algebras and
give an appropriate criterion for the convergence of the associated
spectral sequences. Section \ref{Sec.ACS} is devoted to the
description of the algebras of complete symbols that we are interested
in. In Section \ref{Sec.OA}, we fix the assumptions on our groupoids
and also prove that these algebras depend only on $\pi : M \to B$, as
long as Assumptions \eqref{eq.A1} and \eqref{eq.A2} are
satisfied. Section \ref{Sec.Hochschild} is devoted to the computation
of the Hochschild homology of our algebras of complete symbols. Along
the way, we compute several other homology groups associated to
certain Poissson manifolds. In Section \ref{Sec.Rel}, we extend the
main results of the previous sections to the relative case. The last
section, Section \ref{Sec.Examples}, treats in detail a few
examples. In particular, we obtain a description of the space of
traces on our algebras of complete symbols.  Note that in this paper
almost all results are formulated in the $\ZZ/2\ZZ$-graded case, in
view of some possible applications.

We hope that the results of this paper will find applications to the
index theorem for families \cite{AtiyahSinger4,Bismut} or to its
generalization to families of fibrations by manifolds with boundary
\cite{BismutCheeger}.

We would like to thank A. Connes, R. Lauter, R. Melrose, and
B. Monthubert for useful discussions. The second named author would
like to thank the Max Planck Institute for Mathematics in Bonn for its
kind hospitality while parts of this work were completed. He would
also like to thank M. Crainic for explaining his results with
R. Fernandes on the integration of Lie algebroids
\cite{CrainicFernandes} (see also Section~\ref{Sec.Examples}).

\section{Topological filtered algebras\label{Sec.filtered}}

Topologically filtered algebras were introduced in
\cite{BenameurNistor} to provide a natural framework for the algebras
of complete symbols associated to algebras of pseudodifferential
operators. In this section we review the definition of topologically
filtered algebras and a few other relevant facts. The complexes
computing the various homologies of these algebras have to be defined
appropriately. In view of the applications that we have in mind, we
have found it necessary to extend our setting to include that of
$\ZZ/2\ZZ$-graded algebras.  For basic facts about pseudodifferential
operators, see one of the many nice monographs available
\cite{NazaShatSter}, \cite{Shubin}, or \cite{Taylor}.

We begin by recalling the definitions of Hochschild and cyclic
homology groups of a topological algebra $\alge$. A good reference is
Connes' book \cite{ConnesBOOK}. See also \cite{Karoubi,Loday}.
See \cite{Kassel} for the homology of $\ZZ/2\ZZ$-graded algebras.
These definitions have to be (slightly) modified when the multiplication of
our algebra is only separately continuous. We thus discuss also the
changes necessary to handle the class of algebras that we are interested
in, that of ``topologically filtered algebras'' (Definition
\ref{def.t.filtered}), and then we prove some results on the homology
of these algebras.

First we consider the case of a topological algebra $\alge$. Here
``topological algebra'' has the usual meaning, that is, $\alge$ is a
real or complex algebra, which is at the same time a locally convex
space such that the multiplication $\alge \times \alge \to \alge$ is
continuous when $\alge \times \alge$ is endowed with the product
topology.  Denote by $\potimes$ the projective tensor product and
$\mathcal H_n(\alge):= \alge^{\potimes n +1}$, the completion of
$\alge^{\otimes n+1}$ in the topology of the projective tensor
product. Also, we denote as usual by $\pa a \in \ZZ/2\ZZ$ the degree
of an element in the $\ZZ/2\ZZ$-graded algebra and by $b'$ and $b$ the
Hochschild differentials:
\begin{equation}\label{eq.def.bb'}
\begin{gathered}
 b'(\Tt)=\sum_{i=0}^{n-1} (-1)^ia_0\otimes\ldots\otimes a_i
 a_{i+1}\otimes\ldots\otimes a_n,\\ b(\Tt)=b'(\Tt) +
 (-1)^{n + \mu} a_na_0\otimes\ldots\otimes a_{n-1},
\end{gathered}
\end{equation}
where $\mu = \pa a_n (\pa a_0 + \ldots + \pa a_{n-1})$.

The {\em Hochschild homology groups} of the algebra
$\alge$, denoted $\Hd_*(\alge)$, are then the homology groups of the
complex $({\mathcal H}_n(\alge), b).$ By contrast, the complex
$({\mathcal H}_n(\alge), b')$ is often acyclic, for example when
$\alge$ has a unit. A topological algebra $\alge$ for which
$({\mathcal H}_n(\alge), b')$ is acyclic is called {\em $H$-unital}
(or, better, {\em topological $H$-unital}), following Wodzicki
\cite{Wodzicki}.

We now define cyclic homology. Assume first that $\alge$ is unital. We
shall use the notation of \cite{ConnesNCG}. See also \cite{Karoubi}.

\begin{equation}
\begin{gathered}
s(\Tt)=1\otimes \Tt, \\ t(\Tt)=(-1)^{n + \mu}  a_n\otimes
a_0\otimes\ldots\otimes a_{n-1}, \\ B_0(\Tt)=s\sum_{k=0}^{n}
t^k(\Tt), \;\; \text{ and } \; B = (1 - t)B_0,
\end{gathered}
\end{equation}
where $\mu = \pa a_n (\pa a_0 + \ldots + \pa a_{n-1})$, as above.
Then $[b,B]_+ := bB + Bb = B^2 = b^2 = 0$, and hence, if we define
\begin{equation}
{\mathcal C}(\alge)_n=\bigoplus_{k\geq 0} {\mathcal H}_{n -
2k}(\alge),
\end{equation}
$({\mathcal C}(\alge),b+B),$ is a complex, called {\em the cyclic
complex} of $\alge$, whose homology is by definition the {\em cyclic
homology} of $\alge$, as introduced in \cite{ConnesNCG} and
\cite{Tsygan}. The cyclic homology groups of the algebra $\alge$ are
denoted $\Hc_n(\alge)$. For algebras without unit $\alge$ one
considers the algebra with adjoint unit $\alge^+$ and then the cyclic
homology of $\alge$ is the kernel of the map $\Hc_n(\alge^+) \to
\Hc_n(\CC)$ induced by the augmentation morphism $\alge^+ \to \CC$.

Consideration of the natural periodicity morphism $\mathcal C_n(\alge)
\to \mathcal C_{n-2}(\alge)$ easily shows that cyclic and Hochschild
homology are related by a long exact sequence
\begin{equation} \label{eq.SBI}
\ldots\rightarrow \Hd_{n}(\alge)\stackrel{I}{\longrightarrow}
\Hc_{n}(\alge)\stackrel{S}{\longrightarrow}
\Hc_{n-2}(\alge)\stackrel{B}{\longrightarrow}
\Hd_{n-1}(\alge)\stackrel{I}{\longrightarrow}\ldots\, ,
\end{equation}
with the maps $I$, $B$, and $S$ explicitly determined. The map $S$ is
also called the {\em periodicity} operator.  See
\cite{ConnesNCG,Loday-Quillen1}. This exact sequence exists whether or
not $\alge$ is endowed with a topology.

For the algebras that we are interested in, however, the
multiplication is usually only separately continuous, but there will
usually exist an increasing multi-filtration $F_{n,l}^m \alge \subset
\alge$ of $\alge$,
\begin{equation*}
 F_{n,l}^m \alge \subset F_{n',l'}^{m'}\alge, \quad \text{if }
 n \le n',\, l\le l',\, \text{ and } m \le m',
\end{equation*}
by closed, complemented subspaces, invariant under the $\Z_2$-grading,
and satisfying the following properties:
\begin{enumerate}
\item\ $\alge = \displaystyle{ \cup_{n,l,m} } F^m_{n,l} \alge$;
\item\ The union $\alge_n := \displaystyle{ \cup_{m,l} F_{n,l}^m \alge
}$ is a closed subspace such that
\begin{equation*}
 	F_{n,l}^m \alge = \alge_n \cap \big ( \cup_{j} F_{j,l}^m \alge
 	\big );
\end{equation*}
\item\ Multiplication maps $F_{n,l}^m \alge \otimes F_{n',l'}^m \alge$
to $F_{n + n',l+l'}^m \alge$;
\item\ The maps
$$
 	F_{n,l}^m\alge /F_{n-1,l}^m \alge \otimes F_{n',l'}^m \alge
 	/F_{n'-1,l'}^m \alge \to F_{n+n',l+l'}^m\alge /F_{n + n' -1,
 	l+l'}^m \alge
$$
induced by multiplication are continuous;
\item\ The quotient $F_{n,l}^m\alge/ F_{n-1,l}^m\alge$ is a nuclear  Frechet
space in the induced topology;
\item\ The natural map
$$
 	F_{n,l}^m \alge \to \displaystyle{\lim_{\leftarrow}}\,
 	F_{n,l}^m / F_{n-j,l}^m \alge, \quad j \to \infty
$$
is a homeomorphism; and
\item\ The topology on $\alge$ is the strict inductive limit of the
subspaces $F_{n,n}^n\alge$, as $n \to \infty$ (recall that
$F_{n,n}^n\alge$ is assumed to be closed in $F_{n+1,n+1}^{n+1}\alge$).
\end{enumerate}

\begin{definition}\label{def.t.filtered}\
An algebra $\alge$ satisfying the conditions 1--7, above, will be
called a {\em topologically filtered} algebra.
\end{definition}

It follows from the definition that if $\alge^{[m]} := \displaystyle
{ \cup_{n,l} F^m_{n,l} \alge }$, then $\alge^{[m]}$ is actually a
subalgebra of $\alge$ which is topologically filtered in its own, but
with multi-filtration independent of $m$.

For topologically filtered algebras, the multiplication is not
necessarily continuous, and the definition of the Hochschild and
cyclic homologies using the projective tensor product of the algebra
$\alge$ with itself, as above, is not very useful. For this reason, we
change the definition of the space $\mathcal H_m(\alge)$ to be an
inductive limit:
\begin{equation*}
 \mathcal H_q(\alge) \simeq \displaystyle{\lim_\to} (F_{n,n}^n
 \alge)^{\potimes q+1},
\end{equation*}
the tensor product being the (complete) projective tensor
product.
The Hochschild homology of $\alge$ is then still the homology of the
complex $(\mathcal H(\alge), b)$. Since the projective tensor product
is compatible with the projective limits, we also have
\begin{equation}
 \mathcal H_q(\alge) = \lim_{\to} \left( \lim_{\leftarrow}\,
 (F_{n,n}^n \alge / F_{k,n}^n \alge)^{\potimes q+1} \right ),
\end{equation}
with the induced topology, where first $k \to -\infty$ (in the
projective limit) and then $n \to \infty$ (in the inductive limit).
The operator $B$ extends to a well defined map $B : \mathcal
H_q(\alge) \to \mathcal H_{q + 1}(\alge)$, which allows us to define
the cyclic complex and the cyclic homology of the algebra $\alge$ as
the homology of the complex $(\mathcal C_*(\alge),b+B)$, with
$\mathcal C_q(\alge):= \oplus \mathcal H_{q - 2k}(\alge)$, as for
topological algebras.

We also observe that both the Hochschild and cyclic complexes have
natural filtrations given by
\begin{equation}\label{eq.def.Hfiltr}
 	F_p\mathcal H_q(\alge) := \displaystyle{\lim_{\to}} \left
 	( \displaystyle{\lim_{\leftarrow}} \, \widehat{\otimes}_{j=0}^q
 	\, \big( F_{k_j,m}^m \alge / F_{k,m}^m \alge \big ) \right ),
\end{equation}
where $k_0 + \ldots + k_n \leq p$ defines the filtration. The
projective and inductive limits are such that first $k,l \to - \infty$
(in the projective limit) and then $m \to \infty$ (in the inductive
limit). It follows from the definitions that the operators $S, B$ and
$I$ preserve these filtrations.

For any topologically filtered algebra, we denote
\begin{equation*}
	Gr(\alge) := \oplus_n \, \alge_n/\alge_{n-1}
\end{equation*}
the {\em graded algebra} associated to $\alge$, where $\alge_n$ is the
union $\cup_{l,m}F_{n,l}^m \alge$, as before. Its topology is that of
an inductive limit of Frechet spaces:
\begin{equation*}
	Gr(\alge) \simeq \displaystyle{\lim_{N,m,l \to \infty}}
	\oplus_{n = -N}^{N} F_{n,l}^m\alge/F_{n-1,l}^m\alge.
\end{equation*}
For the algebras like $Gr(\alge)$, we need yet a third way of
topologizing its iterated tensor products. The correct definition is
then
\begin{equation*}
	\mathcal H_q(Gr(\alge)) \simeq \displaystyle{\lim_{N, m ,l \to
	\infty}} (\oplus_{n = -N}^{N}
	F_{n,l}^m\alge/F_{n-1,l}^m\alge)^{\potimes q+1}.
\end{equation*}
The Hochschild homology of $Gr(\alge)$ is the homology of the complex
$(\mathcal H_*(Gr(\alge))$,$b)$. The operator $B$ again extends to a
map $B : \mathcal H_q(Gr(\alge)) \to \mathcal H_{q +
1}(Gr(\alge))$ and we can define the cyclic homology of $Gr(\alge)$ as
above. The operators $S, B$ and $I$ associated to $\mathcal
H_q(Gr(\alge))$ are the graded operators associated with the
corresponding operators
(also denoted $S,B$ and $I$) on  $\mathcal H_q(\alge)$.

The Hochschild and cyclic complexes of the algebra $Gr(\alge)$
decompose naturally as direct sums of complexes indexed by $p \in
\ZZ$. For example, $\mathcal H_q(Gr(\alge))$ is the direct sum of the
subspaces $\mathcal H_q(Gr(\alge))_p$, where
\begin{multline*}
\mathcal{H}_{q}(Gr(\alge))_p = \lim_{m, N ,l \to \infty}
\bigoplus_{k_j} \, \left ( \widehat{\otimes}_{j=0}^n
F^m_{k_j,l}\alge/F^m_{k_j-1,l}\alge \right ), \\ \text{ where
}\; k_0 + k_1 + \ldots + k_n = p \, \text{ and } -N \le k_j
\le N,
\end{multline*}
with the inductive limit topology. The corresponding subcomplexes of
the cyclic complex are defined similarly. We denote by
$\HH_*(Gr(\alge))_p$ and $\HC_*(Gr(\alge))_p$ the homologies of the
corresponding complexes (Hochschild and, respectively, cyclic).

The following two results are well known consequences of standard
results in homological algebra (for topologically filtered algebras
they were proved in \cite{BenameurNistor}).

\begin{lemma}\label{lemma.sp.sq}\ Let $\alge$ be a topologically
filtered algebra.  Then the natural filtrations on the Hochschild and
cyclic complexes of $\alge$ define spectral sequences $\EH^r_{k,h}$
and $\EC^r_{k,h}$ such that
\begin{equation*}
 	\EH^1_{k,h} \simeq \HH_{k+h}(Gr(\alge))_{k} \;\; \text{and }
 	\EC^1_{k,h} \simeq \HC_{k+h}(Gr(\alge))_{k}.
\end{equation*}
Moreover, the periodicity morphism $S$ induces a morphism $S' :
\EC^r_{k,h} \to \EC^r_{k,h-2}$ of spectral sequences. For $r = 1$, the
morphism $S'$ is the graded map associated to the periodicity operator
$S:\HC_{n}(Gr(\alge)) \to \HC_{n-2}(Gr(\alge))$ and the natural
filtration of the groups $\HC_{n}(Gr(\alge))$.
\end{lemma}

\begin{proof}\
The filtration $F_p\maH_q(\alge)$ of the complex computing the
Hochschild homology of $\alge$ gives rise to a spectral sequence
$(E^r)_{r\geq 1}$ with 
\begin{equation*}
	E^1_{k,h} =
	\cohom_{k+h}(F_k\maH(\alge)/F_{k-1}\maH(\alge)),
\end{equation*}
by standard homological algebra. By the definition of the Hochschild
complex of $Gr(\alge)$, we have:
$$
        \cohom_{k+h}(F_k\maH(\alge)/F_{k-1}\maH(\alge)) \simeq
	\HH_{k+h}(Gr(\alge))_k.
$$
This completes the proof for Hochschild homology. The proof for cyclic
homology is similar.
\end{proof}

In our considerations below, we shall need the following classical
result, which was proved for topologically filtered algebras in
\cite{BenameurNistor}. Due to the importance of this result in what
follows and for the convenience of the reader, we include a proof of
it.

\begin{theorem}\label{theorem.conv2}\
Fix an integer $N$ and $a \ge 1$. Let $\alge$ be a topologically
filtered algebra and let $\alge^{[m]} := \cup_{n,l}F_{n,l}^{m}\alge$
be such that $\EH^{a}_{k,h}(\alge^{[m]}) = 0$, for all $k < N$. Then
the spectral sequence $ \EH^r_{k,h} = \EH^r_{k,h}(\alge)$ defined in
Lemma \ref{lemma.sp.sq} converges to $\HH_{k+h}(\alge)$. More
precisely, we have
\begin{equation*}
        \HH_{j}(\alge) \simeq \oplus_{k = N}^\infty \EH^{\infty}_{k,j -
        k}.
\end{equation*}
A similar result holds for the cyclic homology spectral sequence.
\end{theorem}

\begin{proof}\ 
Assume first that the filtration of $\alge$ is independent of $m$. We
have
$$
        \maH(\alge)=\lim_{\leftarrow} \maH(\alge)/F_p\maH(\alge).
$$
This enables to write, for every fixed $q$, the well known associated
$\lim^1$ exact sequence (see \cite{BenameurNistor}[Lemma 6], for
example)
$$
        0\to {\lim_{\leftarrow}}^1
        \cohom_{q+1}(\maH(\alge)/F_p\maH(\alge)) \to
        \cohom_{q}(\maH(\alge)) \to \lim_{\leftarrow}
        \cohom_{q}(\maH(\alge)/F_p\maH(\alge)) \to 0.
$$
Now the spectral sequence $E^r_{k,h}(p)$ associated with the
Hochschild homology of the filtered complex
$\maH(\alge)/F_p\maH(\alge)$ converges because it is the shift of a
first quadrant spectral sequence and we have:
$$
        \cohom_q(\maH(\alge)/F_p\maH(\alge)) \cong \oplus_{t+s=q}
	E^{\infty}_{t,s}(p).
$$
Therefore, the homology groups $\cohom_q(\maH(\alge)/F_p\maH(\alge))$
are endowed with a filtration $\maF_t(p)$ so that
\begin{equation}\label{eq.Ft}
        \maF_t(p)/\maF_{t-1}(p) \cong E^{\infty}_{t,q-t}(p).
\end{equation}
Furthermore, the spectral sequence $E^r_{k,h}(p)$ satisfies:
%\begin{equation}\label{eq.er}
%        E^r_{k,h}(p) = 0 \text{ if } k\leq p \text{ and } E^r_{k,h}(p)
%	= \EH^r_{k,h} \text{ if } k > p+r.
%\end{equation}
\begin{equation}\label{eq.er}
       E^r_{k,h}(p) =  \begin{cases} 0  & 
	\text{ if } k\leq p \; \text{ and } \\
	\EH^r_{k,h} & \text{ if } k > p+r
	\end{cases} 
\end{equation}
Let us choose a particular projective system by setting
$$
        A_n := \cohom_q( \maH(\alge)/F_{N-na}\maH(\alge) ), \; B_n :=
	\maF_{N-na+a}(N-na)\,, \text{ and }\; C_n:=A_n/B_n.
$$
Then the ker-coker lemma \cite{AtiyahMacDonald} for the short exact sequence
$$
        0\to \Pi B_n \hookrightarrow \Pi A_n \to \Pi C_n \to 0
$$
gives rise to the following exact sequence \cite{BenameurNistor}[Lemma
7]:
$$
        0\to \lim_{\leftarrow} B_n \to \lim_{\leftarrow} A_n \to
	\lim_{\leftarrow} C_n \to {\lim_{\leftarrow}}^1 B_n \to
	{\lim_{\leftarrow}}^1 A_n \to {\lim_{\leftarrow}}^1 C_n \to 0.
$$

By conditions \eqref{eq.Ft} and \eqref{eq.er}, the natural map
$A_{n+1} \to A_n$ induces the zero map $B_{n+1} \to B_n$ and an
isomorphism $C_{n+1} \to C_n$, for $n\geq 2$. Therefore we get:
$$
        {\lim_{\leftarrow}}^1 A_n =0 \text{ and } \lim_{\leftarrow}
	A_n = C_{n_0}, \quad \forall n_0\geq 2.
$$
And hence, finally,
$$
        \HH_q(\alge) \cong C_{n_0} = \oplus_{l\geq N}
	E^{\infty}_{l,q-l} \cong \oplus_{l\in \Z}
	\EH^{\infty}_{l,q-l}.
$$
The proof in the general case when the filtration $F_{n,l}^m\alge$ of
$\alge$ depends on $m$ is deduced using that homology is compatible
with taking direct limits.
\end{proof}

The above isomorphism is not natural, in general, but comes from a
filtration $F_k\HH_q(\alge)$ of $\HH_q(\alge)$ whose subquotients
$F_k\HH_q(\alge)/F_{k-1}\HH_q(\alge)$ identify naturally with
$\EH^{\infty}_{k,q-k}$, see \cite{MacLane,MacLaneMoerdijk}.

It is useful to mention here that the composite map
$$
        \Hd_q(\alge) \stackrel{\II}{\rightarrow} \Hc_q(\alge)
        \stackrel{\B}{\rightarrow}\Hd_{q+1}(\alge)
$$
preserves the filtrations and hence it induces natural maps 
\begin{equation*}
	(\B\circ \II)^{(r)} : \EH^r_{k,h}(\alge) \to \EH^r_{k,h +1}(\alge).
\end{equation*}
For $r =1$, this map is the composition of the corresponding morphisms
$$
	\Hd_q(Gr(\alge)) \rightarrow \Hc_q(Gr(\alge))\rightarrow
	\Hd_{q+1}(Gr(\alge))
$$ 
for the graded algebra of $\alge$.

%\section{Families of manifolds with corners  \label{Sec.OA}}
\section{Algebras of complete symbols  \label{Sec.ACS}}

We now introduce the algebras of complete symbols that we study in
this paper.

We shall follow the standard notation for groupoids and Lie
algebroids, using the conventions of \cite{LauterNistor}.  In
particular, if $\GR$ is a differentiable groupoid with space of units
$M$, then $d, r : \GR \to M$ denote the domain and range maps,
respectively, so that the composition $gg'$ of two elements $g,g' \in
\GR$ is defined if, and only if, $d(g) = r(g')$.

We shall also follow \cite{BenameurNistor} for some specific
constructions involving manifolds with corners, some of which are
recalled below. As in that paper, we are interested in certain
specific groupoid algebras associated to manifolds with corners. If
$\GR$ is a differentiable groupoid with space of units $M$ and $E \to
M$ is a $\ZZ/2\ZZ$-graded vector bundle, then we shall denote by
$$
 	\PS{\infty} = \cup_{m\in \ZZ} \PS{m}
$$
the algebra of pseudodifferential operators on $\GR$ acting on
sections of the vector bundle $r^*E$. We also define
$$
 	\PS{-\infty} := \cap_{m\in \ZZ} \PS{m}.
$$
(see \cite{LauterNistor} or \cite{NWX} for definitions).
These two algebras are naturally $\ZZ/2\ZZ$-graded.

We shall denote by $\mO$ the space of smooth functions on the interior
of $M$ that have only Laurent singularities at the boundary faces. If every
hyperface $H$ of $M$ has a defining function $x_H$, then $\mO$ is the
ring generated by $\CI(M)$ and $x_H^{-1}$. Let then
$$
       	\alge(\GR; E) := \PS{\infty}/\PS{-\infty}\quad \text{and}\quad
       	\alge_{\maL}(\GR; E) = \mO \alge(\GR;E).
$$
The $\ZZ/2\ZZ$-grading on $E$ then provides us with a natural
$\ZZ/2\ZZ$-grading on the algebras $\alge(\GR; E)$ and
$\alge_{\maL}(\GR; E)$ too.

\begin{proposition}\label{prop.t.f}\
Assume that $\GR$ and $M$ are as above and that $M$ is
$\sigma$-compact. Then the quotients $\alge(\GR; E)$ and
$\alge_{\maL}(\GR; E)$ are topologically filtered algebras.
\end{proposition}

\begin{proof}\
The proof is exactly the same as the one given for $E = \CC$ in
\cite{BenameurNistor}.
\end{proof}

Let $A(\GR)$ be the Lie algebroid of $\GR$ (see \cite{LauterNistor})
and let $S^*(\GR)$ be the sphere bundle of $A^*(\GR)$, that is, the
set of unit vectors in the dual of the Lie algebroid of $\GR$, and
denote $\cohom_c^{[q]} = \oplus_{k \in \ZZ} \cohom_c^{q + 2k}$.

\begin{theorem}\
Assume that the base $M$ is $\sigma$-compact, then the periodic cyclic
homology of the algebra $\alge(\GR;E)$ is given by
\begin{equation}
 	\Hp_m( \alge(\GR; E) ) \simeq \cohom_c^{[m]}(S^*(\GR) \times
 	S^1).
\end{equation}
\end{theorem}

\begin{proof}\
An argument similar to that of Lemma \ref{Lemma.HH.AL} in Section
\ref{Sec.Hochschild} shows that the Hochschild homology is unchanged
by introducing the extra vector bundle $E$ and the
$\ZZ/2\ZZ$-grading. Standard homological algebra arguments then show
that the same is true for cyclic and periodic cyclic homology. The
result follows then from the case $E = \CC$ that was proved in
\cite{BenameurNistor}.
\end{proof}

To state the result for the algebra $\alge_{\maL}(\GR;E)$, we need
first to recall a construction from \cite{BenameurNistor} that will
be used several times in what follows.

Let $P$ be a manifold with corners. Then $\maL(P)$ is a space
naturally associated to $P$ and defined as follows. Consider for each
face $F$ of $P$ the space $F \times (S^1)^k$, where $k$ is the
codimension of the face. We establish a one-to-one correspondence
between the canonical $k$ copies of the unit circle in $(S^1)^k$ and
the faces $F'$ of $P$ containing $F$, of dimension one higher than
that of $F$. We then identify the points of the disjoint union $\cup F
\times (S^1)^k$ as follows.  If $F \subset F'$ and $F'$ corresponds to
the variable $\theta_i\in S^1$ we identify $(x, \theta_1, \ldots,
\theta_{i-1},1, \theta_{i+1}, \ldots, \theta_{k}) \in F \times
(S^1)^k$ to the point $(x, \theta_1, \ldots, \theta_{i-1},
\theta_{i+1}, \ldots,\theta_{k}) \in F' \times (S^1)^{k-1}$ (same
$x$). The resulting quotient space is by definition $\maL(P)$.

By construction, there exists a continuous map $p_{\maL}:\maL(P) \to
P$. Let $J_\epsilon = S^1 \cup [1, 1+\epsilon) \subset \CC$, for some
$\epsilon >0$, with $S^1$ identified with a subset of the complex
plane.  Then the space $\maL(P)$ is locally modeled by $J_\epsilon^k
\times \RR^{n-k}$, above each point of $P$ belonging to an open face of
codimension $k$. 

Suppose now that $P \to B$ is a fibration by manifolds with corners
with $B$ smooth. In particular, $B$ is a smooth manifold.  Let $Q$ be
the typical fiber of the fibration $P \to B$.  Then we obtain by the
above construction, a
fibration $\maL(P) \to B$, with typical fiber the locally compact
space $\maL(Q)$.

In \cite{BenameurNistor}, the periodic cyclic homology of several
algebras of complete symbols was computed. These results include our
algebras $\alge_{\maL}(\GR; E)$, when $E$ is trivial.
The result is the same in general.

\begin{theorem}\label{periodic}\
For $j=0,1$, we have
$$
        \Hp_j(\alge_{\maL}(\GR; E)) \simeq \cohom_c^{[j]}
        (\maL(S^*(\GR)) \times S^1 ).
$$
\end{theorem}

\begin{proof}\
When $E$ is trivial one-dimensional, this result was proved in
\cite{BenameurNistor}. The general case is proved in the same way,
using the same argument as in the proof of Lemma \ref{Lemma.HH.AL},
which shows that the Hochschild homology of the algebras
$\alge_{\maL}(\GR;E)$ does not depend on the bundle $E$, thanks to the
Morita invariance of Hochschild homology.
\end{proof}

In particular, when $M$ has no corners and $\GR = M \times_B M$, the
algebra $\PS{\infty}$ is the algebra of smooth families of
pseudodifferential operators along the fibers of $\pi : M \to B$ that
have compactly supported Schwartz kernel. Let us denote by $\alge(M
\vert B; E) = \alge(\GR;E)$ the algebra associated to this groupoid .

\begin{corollary}\
For any fibration $\pi : M \to B$ of smooth manifolds (without
boundary), we obtain
$$
	\Hp_j(\alge(M \vert B; E)) \simeq \cohom_c^{[j]}(S^*_{vert}(M)
	\times S^1), \quad j = 0,1.
$$
\end{corollary}

This leads to a complete determination of the periodic cyclic homology
of the algebras $\alge (\GR;E) :=\PS{\infty}/\PS{\infty}$ and
$\alge_{\maL} (\GR;E) := \mO\alge (\GR;E)$. The result is moreover
easily expressed in a uniform manner for all differentiable groupoids
$\GR$. The Hochschild homology of these algebras seems to be more
difficult to compute, in general.  Finding the groups
$\Hd_*(\alge_{\maL}(\GR;E))$ in general seems to depend on finding the
right language in which to express this result.  Needless to say,
finding the right language to express the groups
$\Hd_*(\alge_{\maL}(\GR;E))$ and then determining them in general is a
worthy problem.

We shall determine the groups $\Hd_*(\alge_{\maL}(\GR;E))$ for a class
of groupoids that, roughly speaking, consists of families of groupoids
of the kind considered in \cite{BenameurNistor}. We now proceed to
describe this class in detail.

\section{Families of manifolds with corners  \label{Sec.OA}}

To describe the class of differentiable groupoids $\GR$ for which we
shall determine the groups $\Hd_*(\alge_{\maL}(\GR;E))$, we first
describe the assumptions on the space of units of $\GR$. We shall
denote the space of units of $\GR$ by $M$, where $M$ is a
differentiable manifold, possibly with corners, as before.  We shall
assume that there exists a smooth manifold (without corners) $B$ and a
map $\pi : M \to B$ that makes $M$ a differentiable fiber bundle over
$B$ with fiber $F$.  We regard this fiber bundle as a being the fiber
bundle associated to a principal bundle with structure group
$\Diffeo(F)$, the group of diffeomorphisms of $F$ that map faces to
faces.  {}From now on and throughout the paper, we shall
denote $n = \dim (M)$, $q = \dim (B)$. Also, we shall denote by $p$
the dimension of the fibers of $\pi : M \to B$, so, in particular, $n
= p + q$.

Fix $M$ as above. We now describe our three assumptions on the groupoid
$\GR$. Our {\em first assumption} on $\GR$ is that for any arrow $g \in
\GR$, the domain and range of $g$ are in the same fiber of $\pi : M
\to B$, that is,
\begin{equation}\label{eq.A1}
        \pi( d(g) ) = \pi( r(g) ), \quad \forall g \in \GR.
\end{equation}
The intuitive meaning of this condition is that the natural action of
$\PS{\infty}$ on $\CIc(M)$ via the vector representation
\cite{LMN,NWX} is given by families of operators acting on the fibers
of $\pi$.

Let $T_{vert}M$ be the vertical tangent bundle to the fibration $\pi:
M \to B$. Denote as above by $\mO$ the space of smooth functions on
the interior $M_0$ of $M$ that have only Laurent (or rational) type
singularities at the faces of $M$. Let us denote by $\varrho : A(\GR)
\to TM$ the anchor map of the Lie algebroid of $\GR$. Our {\em second
assumption} on $\GR$ is that the map $\varrho_\Gamma : \Gamma(A(\GR))
\to \Gamma(T_{vert}M)$ defined by $\varrho$ induces an isomorphism
\begin{equation}\label{eq.A2}
        \mO \otimes_{C^{\infty}(M)} \Gamma(A(\GR)) \simeq \mO
        \otimes_{C^{\infty}(M)} \Gamma(T_{vert}M),
\end{equation}
of vector spaces. Clearly the above map preserves the Lie bracket,
so we get an isomorphism of Lie algebras also.

Our {\em third and last assumption} on $\GR$ is a local triviality
condition on the algebra $\alge_{\maL}(\GR; E)$. To state this
assumption, we need to fix the notation. For any open set $V \subset
B$, we denote by $\GR_V$ the reduction of $\GR$ to $\pi^{-1}(V)$. Our
previous assumptions on $\GR$ give that $\GR_V = (\pi \circ
d)^{-1}(V)$.  Similarly, for every point $b \in B$, we denote by
$\GR_b$ the reduction of $\GR$ to $\pi^{-1}(b)$. Again, our
assumptions give us that $\GR_b = (\pi \circ d)^{-1}(b) = (\pi \circ
r)^{-1}(b)$. Let us observe that
$$
        \mathcal O (\pi^{-1}(V)) \otimes \left ( \Psi^{\infty}(V
        \times \GR_b;E)/ \Psi^{-\infty} (V \times \GR_b;E) \right )
$$
has a natural
filtration and a natural completion to a topologically filtered
algebra, denoted
$$
        \mathcal O(\pi^{-1}(V)) \otimes_{tf} \left (\Psi^{\infty}(V
        \times \GR_b;E)/ \Psi^{-\infty} (V \times \GR_b;E) \right).
$$

Then, for any $b \in B$, we assume the existence of an open
neighborhood $V \subset B$ of $b$ and a $\CI(B)$-linear isomorphism
\begin{multline}\label{eq.A3}
        \alge_{\maL}(\pi^{-1}(V)\vert V; E_{\vert V}) := \mathcal O
        (\pi^{-1}(V)) \left ( \Psi^{\infty}(\GR_V;E)/
        \Psi^{-\infty}(\GR_V;E) \right ) \\ \simeq \mathcal
        O(\pi^{-1}(V)) \otimes_{tf}\left ( \Psi^{\infty}(V \times
        \GR_b;E)/ \Psi^{-\infty} (V \times \GR_b;E) \right )
\end{multline}
of topologically filtered algebras, where $V \times \GR_b$ is the
product groupoid, with $V$ consisting of just units and the operations
being defined pointwise.

The three assumptions above, Equations \eqref{eq.A1}, \eqref{eq.A2},
and \eqref{eq.A3} are not independent, as we shall see shortly.  See
\cite{Weinstein} for some basic facts about Poisson manifolds.

\begin{lemma}\label{lemma.Poisson.Iso}\
Assume that \eqref{eq.A2} is satisfied. Then the morphism $\varrho$
above induces an isomorphism
$$
        \mO \otimes_{C^{\infty}(M)} C^{\infty} (A^*(\GR)) \simeq \mO
        \otimes_{C^{\infty}(M)} C^{\infty}(T^*_{vert}M)
$$
of Poisson algebras.
\end{lemma}

\begin{proof}\ Let $X,Y$, and $Z \in \Gamma(A(\GR))$. Then $X,Y$,
and $Z$ define functions (denoted by the same letter) $X, Y, Z :
A^*(\GR) \to \RR$. Assume $Z = [X, Y]$. Then the Poisson bracket on
$C^{\infty} (A^*(\GR))$ is uniquely determined by $\{X,Y\} = Z$. The
equation $\varrho([X,Y]) = [\varrho(X),
\varrho(Y)]$ shows that the natural map 
$$
	C^{\infty}
	(A^*(\GR)) \to C^{\infty}(T^*_{vert}M)
$$ 
is a Poisson map. The proof is completed by including
$\mO$-coefficients.
\end{proof}

Let $M_0:=M \smallsetminus \pa M$ be as above the interior of $M$ and let
$T_{vert}M_0$ be the vertical tangent bundle to the smooth fibration
$M_0 \to B$. Our second assumption implies, in particular, that the
anchor map $\varrho$ restricts to an isomorphism
$$
        A(\GR)\vert_{M_0} \simeq T_{vert}M_0,
$$
of vector bundles.

We now discuss the relation between our three assumptions on $\GR$.
It turns out that these assumptions do not play equal roles. In fact,
the second assumption implies the third one, and, under some weak
assumptions on $\GR$ ($d$-connectivity) it also implies the first
assumption. The following considerations are however somewhat
independent from the rest of the paper, and, for the purpose of
Hochschild homology computations, the reader can just ignore some of
the results below, but instead consider all three assumptions on $\GR$.

First, let us notice that, in the same spirit as the above lemma, we
get an isomorphism of the algebras of differential operators
corresponding to $\GR$ and to $M$. More precisely, let $\Diff(M,\GR)$
be the algebra of differential operators on $M$ generated by $\CI(M)$
and $\Gamma(A(\GR))$. Similarly, let $\Diff(M)$ be the algebra of
differential operators on $M$ generated by $\CI(M)$ and
$\Gamma(TM)$. The anchor map $\varrho$ then gives rise to an inclusion
\begin{equation}
	\varrho_{\Diff} : \Diff(M,\GR) \to \Diff(M).
\end{equation}

\begin{proposition}\label{prop.diff.ops}\
Suppose that the map $\varrho_\Gamma : \Gamma(A(\GR)) \to \Gamma(TM)$
defined by $\varrho$ is injective. Then our second assumption on $\GR$,
Equation \eqref{eq.A2}, is equivalent to the fact that
$\varrho_{\Diff} : \Diff(M,\GR) \to \Diff(M)$ induces an isomorphism
$$\mO\Diff(M,\GR) \to \mO\Diff(M).$$
\end{proposition}

\begin{proof}\ The space of vector fields on a manifold coincides
with the space of first order differential operators without constant
term (\ie that send the function constant equal to 1 to 0). Thus, the
isomorphism $\mO\Diff(M,\GR) \simeq \mO\Diff(M)$ is equivalent to the
fact that $\mO \Gamma(A(\GR))$ maps surjectively onto $\mO
\Gamma(TM)$.  Since this map is injective by assumption, the result
follows.
\end{proof}

The algebras $\alge_{\maL}(\GR;E)$ turn out to depend only on $\pi : M
\to B$.

\begin{theorem}\label{theorem.all.isom}\
The algebras $\alge_{\maL}(\GR;E)$ are independent of $\GR$, as long
as assumptions \eqref{eq.A1} and \eqref{eq.A2} are satisfied.
\end{theorem}

\begin{proof}\ Assume $E$ is trivial, for simplicity. Let 
$\Psi_{prop}^\infty(M\vert B)$ be the algebra of smooth, properly
supported families of operators acting on the fibers of 
\begin{equation*}
	\pi_0 : M_0 := M \smallsetminus \pa M \to B.
\end{equation*} 
Recall \cite{LMN,NWX} that the vector representation $\pi_{v} :
\tPS{\infty} \to
\End(\CIc(M_0))$ is defined uniquely by 
$$
	(\pi_v(P) f) \circ r = P(f \circ r),
$$ 
see \cite{NWX} for details.  Assumption \eqref{eq.A1} shows that
$\pi_v$ factors through a morphism $\tPS{\infty} \to
\Psi_{prop}^\infty(M\vert B)$. Assumption \eqref{eq.A2} implies that
$d^{-1}(x) \simeq \pi_0^{-1}(\pi_0(x))$ for all $x \in M_0$, so the
action of some $P \in \tPS{\infty}$ on $\GR$ is completely determined
by the action of $P$ on the fibers of $\pi_0$. This means that the
vector representation $\pi_v$ of $\tPS{\infty}$ is injective, so
$\tPS{\infty}$ identifies with a subalgebra of
$\Psi_{prop}^\infty(M\vert B)$. Consequently, $\alge_\maL(\GR)$
identifies with a subalgebra of $\mathcal B :=
\Psi_{prop}^\infty(M\vert B)/\Psi_{prop}^{-\infty}(M\vert B)$.

We now argue that Proposition \ref{prop.diff.ops} and asymptotic
completeness implies that the image of $\tPS{\infty}$ in $\mathcal B$
is independent of $\GR$. Indeed, it is enough to check that the image
of $\tPS{m} \to \mathcal B$ is independent of $\GR$, for any $m$. Let
$D \in \mO \Diff(M;\GR)$ be an elliptic differential operator in
$\mO\tPS{k}$, for some fixed $k \ge 1$. Let $Q$ be a parametrix of
$D$. Then Proposition \ref{prop.diff.ops} implies that 
\begin{equation}
	\mO \Diff(M;\GR)[Q] = \mO \Diff(M)[Q].
\end{equation}

Let $\GR_1$ be another differentiable groupoid satisfying Assumptions
\eqref{eq.A1} and \eqref{eq.A2}, then 
\begin{equation*}
	\mO \Diff(M;\GR_1)[Q] = \mO
	\Diff(M;\GR)[Q],
\end{equation*}
also. Because $\mO (\tPS{m}/\tPS{-\infty}) \to \mathcal B$ is
continuous and injective and the image of the space of operators of
order at most $m$ of $\mO \Diff(M;\GR)[Q]$ is dense in
$\mO\tPS{m}/\tPS{-\infty})$, we obtain that the {\em closure} of the
range of $\tPS{m}$ in $\mathcal B$ does not depend on $\GR$. By
looking at the complete symbols of the images of $\mO(\tPS{m})$ and
$\Psi^m(\GR_1)$ in $\mathcal B$ and using asymptotic completeness of
the algebras of pseduodifferential operators $\mO(\tPS{\infty})$ and
$\Psi^\infty(\GR_1)$, we obtain that the actual range of $\tPS{m}$ in
$\mathcal B$ is independent of $\GR$, as desired.
\end{proof}

Let us recall that $\GR$ is {\em $d$-connected} if, and only if, all
the sets $\GR_x :=d^{-1}(x)$ are connected.

\begin{corollary} Suppose $\GR$ is a differentiable groupoid with
units $M$. Then assumption \eqref{eq.A2} implies 
assumption \eqref{eq.A3}. If $\GR$ is also  $d$-connected,
then \eqref{eq.A2} implies also assumption \eqref{eq.A1}.
\end{corollary}

\begin{proof}\ By  Theorem \ref{theorem.all.isom},
it is enough to check \eqref{eq.A3} for any fixed groupoid $\GR$
satisfying \eqref{eq.A1}. In particular, we can choose $\GR$ to be
locally a product, in which case \eqref{eq.A3} is trivially satisfied.
(For example, we could take $\GR = \GR_{M,b}$, the $b$-groupoid
defined in Section \ref{Sec.Examples}.)

Let $X_1, \ldots, X_m$ be sections of $A(\GR)$. We shall write
$\varrho(X_j)$ for $\varrho_{\Gamma}(X_j)$. Then
\begin{equation*}
	\pi(\exp(\varrho(X_1)) \ldots \exp(\varrho(X_m))x) 
	= \pi(x),
\end{equation*}
for all $x \in M$. The assumption that $\GR$ be $d$-connected is
equivalent to the assumption that, for any $g \in \GR$, there exist
$X_1, \ldots, X_m$ as above such that 
$$
	r(g) = \exp(\varrho(X_1)) \ldots \exp(\varrho(X_m)) d(g),
$$
see \cite{Mackenzie}.
\end{proof}

It also follows from the above discussion that it is enough for our
computations to consider a ``typical'' algebra for each fibration $\pi
: M \to B$.  There are several choices of ``typical'' algebras, in
general.  One possible choice, the $b$-calculus, as well as the result
of our computations for these algebras, will be described in Section
\ref{Sec.Examples}.

\section{Hochschild homology of families\label{Sec.Hochschild}}

In this section, we compute the Hochschild homology of the algebras
$\alge_{\maL}(M \vert B ; E)$ introduced in the previous section.
Recall that these algebras are algebras of complete symbols associated
with a groupoid $\GR$ with units $M$ and a fibration $\pi : M \to B$
by manifolds with corners satisfying the assumptions of Equations
\eqref{eq.A1} and \eqref{eq.A2}. The results of this section are
already interesting when the manifold $M$ has no boundary.  Recall
that $n = \dim (M)$, $q = \dim (B)$, and $n = p + q$.

In addition to helping us eliminate our third assumption on $\GR$,
Equation \eqref{eq.A3}, the introduction of the Laurent-type factors
also simplifies the calculations, as in \cite{MelroseNistor} and
\cite{BenameurNistor}. When $B$ is reduced to a point $*$, this also
ensures that the Hochschild homology of $\alge_{\maL}(M) :=
\alge_{\maL}(M\vert *)$ is finite dimensional. For example, the
dimension of the space of traces on $\alge_{\maL}(M\vert *)$ is the
number of minimal faces of $M$ \cite{BenameurNistor}.  Moreover, the
``cone algebras'' described for example in
\cite{SchulzeSterninShatalov} are more closely related to the algebras
$\alge_{\maL}(M \vert B; E)$ than to the algebras $\alge(M \vert B;
E)$. See also \cite{LauterSeiler,Schrohe}.

Our computations will use the Poisson structure of $A^*(\GR)$ and,
more precisely, the ``homogeneous Laurent-Poisson homology'' of
$A^*(\GR) \smallsetminus 0$, where $A^*(\GR) \smallsetminus 0$ is the
dual of the Lie algebroid of $\GR$, with the zero section removed. The
homogeneous Laurent-Poisson homology of $A^*(\GR) \smallsetminus 0$
is defined below and will be identified in terms of the ``homogeneous,
vertical Laurent-de Rham cohomology'' of the fibration $A^*(\GR)
\smallsetminus 0 \to B$ (this cohomology is also defined below).  The
homogeneous Laurent-Poisson homology and the homogeneous vertical
Laurent-de Rham cohomology are natural analogues of the Poisson and,
respectively, de Rham cohomology, which are obtained, roughly
speaking, by introducing Laurent type singularities at the corners of
$M$ and by considering homogeneous forms (on $A^*(\GR) \smallsetminus
0$, for example). See \cite{Monnier,WeinsteinXu,PingXu} for more
on Poisson cohomology.

We begin with the definition of the groups $\cohom_{c,\maL}^{i,j}
(A^*(\GR) \smallsetminus 0 \vert B)_l$, the homogeneous, vertical
Laurent-de Rham cohomology of the fibration $A^*(\GR) \smallsetminus 0
\to B$.   Then we shall discuss Poisson homology and its
variant, the homogeneous Laurent-Poisson homology.

Let $X$ be a manifold with corners and let $\pi_0 : X \to B$ be a
fiber bundle with $B$ smooth.  Let us call the sections of $T_{vert}X$
vertical vector fields, as it is customary. Then the sections of the
dual $T_{vert}^*X$ are called vertical differential forms. There
exists a natural (\ie independent of any choices) differential
operator
\begin{equation}\label{eq.deRham1}
 	d_{vert} : \Gamma(\Lambda^k T_{vert}^*X) \to
 	\Gamma(\Lambda^{k+1} T_{vert}^*X),
\end{equation}
the {\em vertical  de Rham differential}.

Every vertical vector field on $X$ is also a vector field on $X$ in
the usual sense. On the other hand, a form on $X$ restricts to a
vertical form on $X$. Moreover, every vertical form on $X$ is the
restriction of a form on $X$, but we cannot choose that form in a
canonical way. A convenient way to choose extensions of vertical forms
is to consider a splitting of $TX$ into vertical and horizontal
parts. We shall hence fix from now an isomorphism (or splitting)
\begin{equation}\label{eq.splitting}
 	\Theta : TX \simeq T_{vert}X \oplus \pi_0^*TB.
\end{equation}

The splitting $\Theta$ of Equation \eqref{eq.splitting} gives rise to
an embedding $\Theta_k : \Gamma(\Lambda^k T_{vert}^*X) \to
\Omega^k(X)$.  More
generally, we get isomorphisms
$$
	\Lambda^k T^*X \simeq \bigoplus_{i+j=k} \Lambda^i T^*_{vert}X
 	\otimes \pi_0^* \Lambda^j T^*B.
$$
Let $\Omega^{i,j}(X) := \Gamma(X, \Lambda^i T^*_{vert}X \otimes
\pi_0^* \Lambda^j T^*B).$ Then
%\begin{equation*}
 	$\Omega^k(X) \simeq \bigoplus_{i+j=k} \Omega^{i,j}(X),$
%\end{equation*}
and we also have isomorphisms
\begin{equation}\label{eq.tens.bas}
 	\Omega^{i,0}(X) \otimes_{\CI(B)} \Omega^j(B) \ni \omega
 	\otimes \eta \to \omega \wedge \pi_0^*\eta \in \Omega^{i,j}(X).
\end{equation}

The embedding $\Theta_k$ can then be used to define a map $d_{vert}:
\Omega^{i,0}(X) \to \Omega^{i+1,0}(X)$ (using the same notation is
unlikely to cause any confusion in our case). We extend then
$d_{vert}$ to a map
\begin{equation*}
 	d_{vert} : \Omega^{i,j}(X) \to \Omega^{i+1,j}(X),
\end{equation*}
by using the isomorphisms of the Equation \eqref{eq.tens.bas} above
and setting
$$
        d_{vert}( \omega \wedge \pi_0^* \eta ) = d_{vert}(\omega) \wedge
        \pi_0^* \eta
$$
if $\eta \in \Omega^{i,0}(X)$ and $\omega \in
\Omega^j(B)$. Clearly $d_{vert}^2 = 0$. The extension $d_{vert}$ that
we obtain depends on the splitting $\Theta$ of Equation
\eqref{eq.splitting}.  The isomorphism class of the resulting complex,
however, does not depend on $\Theta$.

Let us denote by $\Omega_{\maL}^k(X)$ the space of $k$-differential
forms on the interior of $X$ that have only rational (or Laurent)
singularities near the corners. We shall sometimes call forms with
these properties {\em Laurent}-differential forms. The above
definitions and properties extend to $\Omega_{\maL}^k(X)$ as
follows. Let $\Omega_{\maL}^{i,j}(X) := \mathcal O(X) \Gamma(X,
\Lambda^i T^*_{vert} X \otimes \pi_0^* \Lambda^j T^*B).$ Then
\begin{equation*}
       	\Omega_{\maL}^k(X)\simeq \bigoplus_{i+j=k}
       	\Omega_{\maL}^{i,j}(X),
\end{equation*}
and, as before, we obtain a differential
\begin{equation*}
        d_{vert} : \Omega_{\maL}^{i,j}(X) \to \Omega_{\maL}^{i+1,j}(X).
\end{equation*}
We shall denote by $\cohom_{\maL}^{i,j}(X) = \ker (d_{vert}) /
d_{vert}\Omega_{\maL}^{i-1,j}(X)$ the homology of the above
complex. Similarly, if compactly supported forms are considered, we
obtain a complex whose homology we denote by
$\cohom_{c,\maL}^{i,j}(X)$.

Define the horizontal differential $d_{hor}: \Omega_{\maL}^{i,j}(X)
\to \Omega_{\maL}^{i,j+1}(X)$ as the component of bidegree
$(0,1)$ of $d$.
% CORRECTED
%by:
%$$
%	        d_{hor} (\omega \wedge \pi_0^*\eta) = \omega \wedge
%	        \pi_0^*(d\eta).
%$$
Then $\pa:=d-d_{vert}-d_{hor}$ is known to be a differential and to
have bidegree $(-1,2)$. See for instance \cite{Tondeur}.

The equality  $d^2=0$ is equivalent to the following relations:
\begin{multline*}
	d_{vert}d_{hor}+d_{hor}d_{vert}=0, \quad d_{hor}^2 + \pa
	d_{vert} + d_{vert}\pa =0 \\ d_{vert}^2=0, \quad \text{and }\,
	\pa d_{hor} + d_{hor}\pa =0.
\end{multline*}
The vertical Laurent-de Rham cohomology can be computed in a fairly
explicit manner. Indeed, let $\mathcal F^k$ be the local coefficient
system determined by the Laurent cohomology groups of the fibers of $X
\to B$. Thus $\mathcal F^k$ is a canonically flat vector bundle over
$B$ whose fiber at $b \in B$ is
\begin{equation}\label{maF}
        \mathcal F^k(b) = \cohom_{c,\maL}^k (\pi_0^{-1}(b)) :=
        \cohom_{c,\maL}^k (\pi_0^{-1}(b)\vert b).
\end{equation}
Let $\Omega_c^k(B)$ be the space of compactly supported $k$-forms on
$B$.

\begin{proposition}\label{prop.ind.vert}\
Using the above notation, we have that
\begin{equation}\label{isomo}
        \cohom^{k,h}_{c,\maL}(X \vert B)
        \simeq \Omega_c^h(B) \otimes_{\CI(B)} \Gamma ( \mathcal F^k )
        =: \Omega_c^h(B; \mathcal F^k),
\end{equation}
the space of compactly supported $h$-forms on $B$ with values in
$\mathcal F^k$. In particular, the vertical Laurent-de Rham cohomology
groups $\cohom^{k,h}_{c,\maL}(X \vert B)$ are independent of the
splitting $TX \simeq T_{vert}X \oplus \pi_0^*TB$ used to define them,
see Equation \eqref{eq.splitting}. 

In addition, the action induced by  the horizontal de
Rham differential $d_{hor}$  on
$\cohom^{k,h}_{c,\maL}(X \vert B)$ is isomorphic under \eqref{isomo}
to the de Rham differential on $B$ with coefficients in the locally
constant sheaf $\mathcal F^k$.
\end{proposition}

\begin{proof}\ The above formula is checked right away when
$\pi : X \to B$ is a trivial fiber bundle (\ie $X = B \times F$) by
using the K\"{u}nneth formula for the tensor product of complexes of
nuclear vector spaces, \cite{KaroubiCyclic}. 
%(TO CHECK KAROUBI'S PAPER) 
Moreover, an automorphism of the trivial fiber bundle $X = B
\times F$ does not affect the isomorphism of the proposition. A
partition of unity argument then completes the proof.
\end{proof}

Let $A^*(\GR) \smallsetminus 0$ be obtained from the vector bundle
$A^*(\GR)$, as before, by removing the zero section. We are interested
in the above constructions when $X = A^*(\GR) \smallsetminus 0$ and
$\pi_0 : A^*(\GR) \smallsetminus 0 \to B$ is obtained from the
composition of the maps $A^*(\GR) \to M$ and $\pi : M \to B$.  More
precisely, for us, the relevant cohomology groups are the cohomology
groups obtained by considering homogeneous forms. Let then
$\Omega_{rc,\maL}^{i,j} (A^*(\GR)\smallsetminus 0)_l$ be the space of
$l$-homogeneous forms in
$$
        \Omega_{\maL}^{i,j}(A^*(\GR) \smallsetminus 0) = \mO \Gamma
        (A^*(\GR) \smallsetminus 0,\, \Lambda^i T_{vert} (A^*(\GR)
        \smallsetminus 0)^* \otimes \pi_0^* \Lambda^j T^*B)
$$
whose support projects onto a compact subset of $M$. Here the
homogeneity is considered with respect to the natural action of
$\RR^*_+$ on $A^*(\GR) \smallsetminus 0 $ by dilations. We denote then
by
$$
        \cohom_{c,\maL}^{i,j} (A^*(\GR) \smallsetminus 0 \vert B)_l
$$
the homology of the complex $\Omega_{rc, \maL}^{i,j} (A^*(\GR)
\smallsetminus 0 )_l$ with respect to the vertical de Rham
differential $d_{vert}$.  We shall call these groups the {\em
homogeneous, vertical Laurent-de Rham cohomology groups of}
$A^*(\GR)$.

Similar constructions and definitions are obtained with $T_{vert}^*M$
in place of $A^*(\GR)$. Our second assumption on the groupoid $\GR$,
Equation \eqref{eq.A2}, gives that the two cohomologies are
isomorphic.

\begin{lemma}\label{lemma.red.T}\
The anchor map $\varrho : A(\GR) \to T_{vert}M$ induces a natural
isomorphism
$$
        \cohom_{c,\maL}^{i,j} (A^*(\GR) \smallsetminus 0 \vert B)_l
        \simeq \cohom_{c,\maL}^{i,j} (T_{vert}^*M \smallsetminus 0 \vert
        B)_l.
$$
These groups vanish if $l \not = 0$ and, for $l = 0$, we have 
\begin{multline*}
	\cohom_{c,\maL}^{i,j} (A^*(\GR) \smallsetminus 0 \vert B)_0
       	\simeq \cohom_{c,\maL}^{i,j} (T_{vert}^*M \smallsetminus 0
       	\vert B)_0 \\ \cong \cohom_{c,\maL}^{i,j} (S_{vert}^*(M)
       	\times S^1 \vert B) \cong \cohom_{c,\maL}^{i,j}
       	(S^*(\GR) \times S^1 \vert B).
\end{multline*}
\end{lemma}

\begin{proof}\
The map $\varrho$ induces an isomorphism of the corresponding
complexes, by Equation \eqref{eq.A2} and the definition of the spaces
$\Omega_{rc, \maL}^{i,j} (A^*(\GR) \smallsetminus 0 )_l$.  The
vanishing of the groups $\cohom_{c,\maL}^{i,j} (T_{vert}^*
\smallsetminus 0 \vert B)_l$ for $l \not = 0$ follows from the homotopy 
invariance of de Rham cohomology. The computation of the
$0$-homogeneous cohomology spaces is classical, see for instance
\cite{BrylinskiGetzler,MelroseNistor}.
\end{proof}

We can thus replace $A^*(\GR)$ with $T^*_{vert}M$ for the rest of our
computations of the homogeneous, vertical Laurent-de Rham cohomology
groups of $A^*(\GR)$.

The homogeneous, vertical Laurent-de Rham cohomology can be computed
using a method similar to the one we used to determine the
non-homogeneous homology. Indeed, let $\mathcal F^k$ be the local
coefficient system determined by the Laurent cohomology groups of the
fibers of $\pi_0 : S^*_{vert}(M) \times S^1 \to B$. Thus $\mathcal
F^k$ is a canonically flat vector bundle over $B$ whose fiber at $b
\in B$ is
\begin{equation}\label{maF2}
       	\mathcal F^k(b) = \cohom_{c,\maL}^k (\pi_0^{-1}(b)) =
       	\cohom_{c}^k (\mathcal L(\pi_0^{-1}(b))).
\end{equation}

\begin{proposition}\label{prop.ind.vert2}\
Using the above notation, we have that
\begin{equation}\label{isomo2}
        \cohom^{k,h}_{c,\maL}(A^*(\GR) \smallsetminus 0 \vert B)_0
        \simeq \Omega_c^h(B) \otimes_{\CI(B)} \Gamma (\mathcal F^k) =:
        \Omega_c^h(B; \mathcal F^k).
\end{equation}
For $l \not = 0$, the groups $\cohom^{k,h}_{c,\maL}(A^*(\GR)
\smallsetminus 0 \vert B)_l$ vanish.  The horizontal de Rham
differential $d_{hor}$ induces a differential on
$\cohom^{k,h}_{c,\maL}(X \vert B)$ which is isomorphic under
\eqref{isomo2} to the de Rham differential on $B$ with coefficients in
the locally constant sheaf $\mathcal F^k$.
\end{proposition}

\begin{proof}\ The proof is completely similar to that of Proposition
\ref{prop.ind.vert}. The vanishing of the groups
$\cohom^{k,h}_{c,\maL}(A^*(\GR) \smallsetminus 0 \vert B)_l$ for $l
\not = 0$ follows again from the homotopy invariance of Laurent-de
Rham cohomology.
\end{proof}

Let us now introduce the Poisson homology groups that we are
interested in. The following considerations apply to any regular
Poisson structure. Recall that the Poisson structure on $A^*(\GR)$ is
defined by a two tensor
\begin{equation*}
 	G \in \CI(A^*(\GR), \Lambda^2 T A^*(\GR))
\end{equation*}
so that $\{f,g\}=i_G(df \wedge dg).$ Clearly, the tensor $G$ must
satisfy some non-trivial conditions for the map $\{\, , \, \}$ to
satisfy the Jacobi identity. These conditions turn out to be
equivalent to $[G,G]_{SN} = 0$, where $[\, ,\, ]_{SN}$ is the
Schouten-Nijenhuis bracket \cite{VaismanBook}. The formula for the Poisson
bracket is determined in terms of the Lie algebra structure on the
space of sections of $A(\GR)$. (This was recalled in Lemma
\ref{lemma.Poisson.Iso}.)

Let $i_G : \Omega^k(A^*(\GR)) \to \Omega^{k-2}(A^*(\GR))$ be the
contraction with the tensor $G$. Then we obtain as in \cite{Brylinski}
a differential
\begin{equation}\label{eq.Brylinski}
        \delta := i_G \circ d - d \circ i_G :\Omega^k(A^*(\GR)) \to
        \Omega^{k-1}(A^*(\GR)).
\end{equation}
Explicitly, for any $(f_0, \ldots ,f_k)\in \CI(A^*(\GR))^{k+1}$, the
differential $\delta$ is given by the formula
\begin{multline}
\label{eq.explicit}
        \delta(f_0 d f_1 d f_2 \ldots d f_k) = \sum_{1\leq j\leq k}
        (-1)^{j+1} \{f_0,f_j\}d f_1 \ldots {\widehat{d f_j}} \ldots d
        f_k \\ + \sum_{1\leq i<j\leq k} (-1)^{i+j}f_0d \{f_i,f_j\} d
        f_1 \ldots {\widehat{d f_i}} \ldots {\widehat{d f_j}} \ldots d
        f_k.
\end{multline}

When $\GR$ is our groupoid and under the Assumptions \eqref{eq.A1} and
\eqref{eq.A2}, we can use the splitting $\Theta$ of Equation
\eqref{eq.splitting}. One then shows that $\delta$ decomposes into the
sum of two bihomogeneous differentials \cite{VaismanBook}
$$
	\delta = \delta_{vert} + \alpha,
$$
where $\delta_{vert}:= i_G \circ d_{vert} - d_{vert} \circ i_G$ is the
vertical Poisson differential. The vertical Poisson differential has
bidegree $(-1,0)$. The extra term $\alpha$ has bidegree $(-2,+1)$ and
is in fact given by \cite{BenameurNistorFol}:
$$
	\alpha = i_G \circ d_{hor} - d_{hor} \circ i_G.
$$
See \cite{BenameurNistorFol}. In particular, the commutator $[i_G, \pa]$ is
trivial. It is straightforward to see that the restriction of
$\delta_{vert}$ to vertical differential forms is given by the
following local expression
\begin{multline}
\label{eq.explicit.vert}
	\delta_{vert}(f_0 d_{vert} f_1 d_{vert} f_2 \ldots d_{vert}
	f_k) \\ = \sum_{1\leq j\leq k} (-1)^{j+1} \{f_0,f_j\}d_{vert}
	f_1 \ldots {\widehat{d_{vert} f_j}} \ldots d_{vert} f_k \\ +
	\sum_{1\leq i<j\leq k} (-1)^{i+j}f_0d_{vert} \{f_i,f_j\}
	d_{vert} f_1 \ldots {\widehat{d_{vert} f_i}} \ldots
	{\widehat{d_{vert} f_j}} \ldots d_{vert} f_k.
\end{multline}

 The above formula determines
$\delta_{vert}$ on $\Omega^{j,0}(A^*(\GR))$. To determine
$\delta_{vert}$ in general, we can use the following lemma.

\begin{lemma}\label{lemma.Bilinear}\
Let $\alpha \in \Omega^i(A^*(\GR))$, let $\beta \in \Omega^j(B)$, and
let $\pi_0 : A^*(\GR) \to B$ be the composite projection.
Then
$$
        \delta(\alpha \wedge \pi^*_0(\beta)) =
        \delta(\alpha) \wedge \pi^*_0(\beta)
	\quad \text{and }\; \delta_{vert}(\alpha \wedge \pi^*_0(\beta)) =
        \delta_{vert}(\alpha) \wedge \pi^*_0(\beta).
$$
\end{lemma}

\begin{proof}\ 
It is enough to check the first equation when $\beta = g$ or $\beta =
dg$, for some smooth function $g$ on $B$. Our claim then follows from
the fact that $\{ f, g \circ \pi_0 \} = 0$, 
%(VICTOR, WHY IS THIS
%TRUE--MOULAY: ONE WAY TO PROVE THIS IS TO RELATE IT TO THE COMMUTATOR
%OF PSEUDODIFFERENTIAL FAMILES, BUT THEN $g \circ \pi_0$ COMMUTES WITH
%ANY FAMILY) 
for any smooth function $f$ on $M$ and from the explicit formula for
$\delta$, Equation \eqref{eq.explicit}.

The equation for $\delta_{vert} := i_G \circ d_{vert} - d_{vert} \circ
i_G$ follows from the equation for $\delta$ by checking bidegrees.
\end{proof}

The formula of Equation \eqref{eq.explicit} is valid also when $M$ has
corners and it is easy to check that the differential $\delta$ is
homogeneous of degree $-1$ with respect to the action of $\R^*_+$ on
$A^*(\GR) \smallsetminus 0$. Let $\Omega^k(A^*(\GR) \smallsetminus
0)_l$ be the space of $k$-forms on $A^*(\GR) \smallsetminus 0$ that
are homogeneous of order $l$. We hence obtain a differential
\begin{equation*}
 \delta : \Omega^k(A^*(\GR)\smallsetminus 0)_{l} \to
 \Omega^{k-1} (A^*(\GR)\smallsetminus 0)_{l-1}.
\end{equation*}

Let $\Omega_{rc}^k(A^*(\GR) \smallsetminus 0)_l$ be the subspace of
$\Omega^k(A^*(\GR) \smallsetminus 0)_l$ consisting of forms whose
support projects onto a compact subset of $M$, as before. Because
$\delta$ preserves the support, it maps the space
$\Omega_{rc}^k(A^*(\GR)\smallsetminus 0)_{l}$ to
$\Omega_{rc}^{k-1}(A^*(\GR)\smallsetminus 0)_{l-1}.$

The same result holds with $\delta_{vert}$ and we have:
$$
	\delta_{vert} : \Omega_{{rc,\maL}}^k(A^*(\GR) \smallsetminus
	 0)_l \to \Omega_{{rc,\maL}}^{k-1}(A^*(\GR)\smallsetminus
	 0)_{l-1}.
$$
We obtain in this way a direct sum of complexes $(\maP^k)_{k\in \Z}$
\begin{equation}\label{delta}
 	\maP^k:\quad 0 \to {\mathcal P}^{k}_{2p+q}
 	\stackrel{\delta}{\rightarrow} {\mathcal P}^{k}_{2p+q -
 	1}\stackrel{\delta} {\rightarrow} \ldots
 	\stackrel{\delta}{\rightarrow} {\mathcal P}^{k}_{-k} \to 0,
\end{equation}
where ${\mathcal P}^{k}_{l} = \Omega_{{rc,\maL}}^{k+l}(A^*(\GR)
\smallsetminus 0)_l$.

We shall denote the homology groups of the above complex by
\begin{equation*}
 	\cohom_{\maL,k+l}^{\delta}(A^*(\GR)\smallsetminus 0\vert B)_l
 	:= {\frac{\ker(\delta:{\mathcal P}^{k}_{l}\to {\mathcal
 	P}^{k}_{l-1})} {\delta(\mathcal P^{k}_{l+1})}}.
\end{equation*}
 In the same way we define the vertical  
{\em homogeneous Laurent-Poisson homology} groups using
$\delta_{vert}$ instead of $\delta$ and denote them by 
$$
	\cohom_{\maL,k+l}^{\delta_{vert}}(A^*(\GR)\smallsetminus
	0\vert B)_l:= {\frac{\ker(\delta_{vert}:{\mathcal
	P}^{k}_{l}\to {\mathcal P}^{k}_{l-1})} {\delta_{vert}(\mathcal
	P^{k}_{l+1})}}.
$$
Furthermore, we define for any $(i,j)$
$$
	\maP_l^{i,j} :=
	\Omega^{i+l,j}_{rc,\maL}(T^*_{vert}M\smallsetminus 0)_l.
$$
{}From the results of Section \ref{Sec.OA}, we deduce that
$$
	\maP^k_l \simeq  \bigoplus_{i+j=k} \maP^{i,j}_l.
$$
Note that with respect to the splitting \eqref{eq.splitting}, we have:
$$
	\delta_{vert}:\maP^{i,j}_l  \to \maP^{i,j}_{l-1},
$$
and the vertical Laurent Poisson homology can be computed by fixing
$(i,j)$ and restricting to each $\maP^{i,j}:=\oplus_l
\maP^{i,j}_l$. However, the extra term $\alpha$ does not preserve
$\maP^{i,j}$, and sends $\maP^{i,j}$ to $\maP^{i-1,j+1}$.

The relevance for us of Poisson homology, in general, and of
homogeneous Laurent-Poisson homology, in particular, is that they are
related to the Hochschild homology groups of the algebras
$\alge_{\maL}(M \vert B; E) := \mO (\Psi^{\infty}(\GR;E) /
\Psi^{-\infty}(\GR;E))$ introduced in the previous section, where $E$
is a $\ZZ/2\ZZ$-graded vector bundle.  The following lemma makes this
connection precise.

\begin{lemma}\label{Lemma.HH.AL}\
The algebra $\alge_{\maL}(M \vert B; E)$ is topologically $H$-unital.
The $\EH^2$ term of the spectral sequence associated to
$\alge_{\maL}(M \vert B; E)$ by Lemma \ref{lemma.sp.sq} is given by
\begin{equation*}
 	\EH_{k,h}^2 \simeq \cohom^\delta_{{\maL},k+h}( A^*(\GR)
 	\smallsetminus 0 \vert B)_k.
\end{equation*}
\end{lemma}

\begin{proof}\
When $E$ is a trivial one-dimensional vector bundle and $\GR$ is an
arbitrary groupoid, the above proposition was proved in
\cite[Proposition 7]{BenameurNistor}.

The extension to a non-trivial vector bundle and the $\ZZ/2\ZZ$-graded
case is obtained as a consequence of K\"{u}nneth formula as follows.
Let us assume first that $E$ is trivially graded. The graded algebra
of $\alge_{\maL}(M \vert B; E)$ is the algebra generated by the
homogeneous sections of the lift of $E$ to $A^*(\GR) \smallsetminus
0$. All these algebras are Morita equivalent, so their Hochschild
homology is the same as that of the algebra corresponding to  a
trivial line bundle. Since the Hochschild homology of a commutative
algebra can be represented as the homology of the space of global
sections of a complex of sheaves, the statement about Morita
equivalence is a local statement, which can then be proved locally
using the K\"{u}nneth formula. The spectral sequence of Lemma
\ref{lemma.sp.sq} then tells us that the Hochschild homology of the
algebra $\alge_{\maL}(M \vert B; E)$ is independent of $E$.

In general, let $E = E_+ \oplus E_-$ be the decomposition of $E$ into
the direct sum of the $+1$ and, respectively, $-1$ eigenvalue of the
grading automorphism. As above, we observe that the statement of the
lemma is local, so we can assume that $E = E_+ \oplus E_-$ is such
that both $E_{+}$ and $E_{-}$ are trivial bundles.  Denote by $N$ the
rank of $E$. Then $\alge_{\maL}(M \vert B; E) \simeq M_N(\alge_{\maL}(M
\vert B)).$
%%% Careful here, no bundle $E$

Let $A$ be a (topologically filtered) algebra $A$ and $N$ an
integer. Assume that the grading automorphism of the algebra $M_N(A)$
is given by conjugation with a matrix in $M_N(\CC)$. Then $\Hd_*(M_N(A))
\simeq \Hd_*(A)$.  This follows for example from the K\"{u}nneth formula
in ($\ZZ/2\ZZ$-graded) Hochschild homology (see \cite{Kassel}).
\end{proof}

As with the homogeneous, vertical de Rham homology, we can replace
$A^*(\GR)$ in $\cohom^\delta_{{\maL},k}(A^*(\GR) \smallsetminus 0
\vert B)_l$ with $T_{vert}^*M$, the vertical cotangent bundle to $\pi
: M \to B$.

\begin{lemma}\label{lemma.is.po}\
The anchor map $\varrho : A(\GR) \to T_{vert}M$ induces an isomorphism
$$
       	\cohom^\delta_{{\maL},k} ( A^*(\GR) \smallsetminus 0 \vert
       	B)_l \simeq \cohom^\delta_{{\maL},k} ( T_{vert}^*M
       	\smallsetminus 0 \vert B)_l.
$$
\end{lemma}

\begin{proof}\
This follows from Lemma \ref{lemma.Poisson.Iso} and the explicit
formula for the Poisson bracket, Equation \eqref{eq.explicit}.
\end{proof}

We now proceed as usual and construct a chain map $*_{vert}$ from the
complex that defines vertical Poisson homology to the complex that defines de
Rham cohomology. The chain map $*_{vert}$ is, in a certain sense, a
vertical symplectic $*$-operator. It corresponds to the canonical
symplectic forms on the cotangent spaces of the fibers of $\pi : M \to
B$. Denote by $M_b := \pi^{-1}(b)$, $b \in B$, and by $\omega_b$ the
symplectic form on $T^*M_b$. There exists then a 2-form $\omega$ on
$T_{vert}^*M : = \cup_{b \in B} T^*M_b$ that restricts on each fiber
$T^*M_b$ to the form $\omega_b$. This form is certainly not unique.
There will be, however, a unique form $\omega \in
\Omega_{\maL}^{2,0}(T_{vert}^*M)$ with this property.  We shall call
this form $\omega$ the {\em vertical symplectic form} of
$T_{vert}^*M$. It depends on the splitting of Equation
\eqref{eq.splitting}.

The vertical symplectic volume form on $T^*_{vert}M$ is defined by
analogy to be $vol_{vert}(M):= \omega^{p}/p!$. Next, we define
$*_{vert} : \Omega_{\maL}^{k,0}(T^*_{vert}M) \to
\Omega_{\maL}^{2p-k,0}(T^*_{vert}M)$ by the equation
$$
        \beta\wedge (*_{vert}\alpha) = (\beta,\alpha)_\omega \cdot
        vol_{vert}(M), \quad \forall \, \alpha,\beta \in
        \Omega_{\maL}^{k,0}(T^*_{vert}M),
$$
where $(\;,\;)_\omega$ is the bilinear form induced by the symplectic
form. Then we obtain that
$*_{vert}(\Omega_{\maL}^{k,0}(T^*_{vert}M)_l) =
\Omega_{\maL}^{2p-k,0}(T^*_{vert}M)_{l+p-k}$. Finally, to define
\begin{equation*}
        *_{vert} : \Omega_{\maL}^{i,j}(T^*_{vert}M) \to
         \Omega_{\maL}^{2p-i,j}(T^*_{vert}M)
\end{equation*}
in general, it is enough to define $*_{vert}\alpha$ when $\alpha =
\eta \wedge \pi_0^*\beta$, with $\eta \in \Omega_{\maL}^{i,0}(M)$ and
$\beta \in \Omega^j(B)$, where $\pi_0 : T^*_{vert}M \to B$ is the
induced projection. We set then
\begin{equation}\label{def.*}
        *_{vert}(\alpha)= *_{vert}(\eta \wedge \pi_0^*\beta):=
        *_{vert}(\eta) \wedge \pi_0^*\beta.
\end{equation}
Similarly, we obtain again that
$*_{vert}(\Omega_{\maL}^{i,j}(T^*_{vert}M)_l) =
\Omega_{\maL}^{2p-i,j}(T^*_{vert}M)_{l+p-i}$.

The usual properties of the symplectic $*$-operator in relation to de
Rham and Poisson homology extend to $*_{vert}$.

\begin{proposition}\ Let $*_{vert}$ be the operator defined above. Then

(1)\ $*_{vert}^2=id$;

(2)\ $(-1)^{i+1} d_{vert} \circ *_{vert} = *_{vert} \circ
\delta_{vert}$ on $\Omega^{i,j}(T^*_{vert}M)$.\\
We can extend the range of both formulas to include homogeneous forms
or forms with Laurent type singularities.
\end{proposition}

\begin{proof}\
Both formulas are well known when $B$ is reduced to a point
\cite{Brylinski}. Using a version with parameters of this particular
case, we obtain that the two formulas are correct on
$\Omega^{i,0}(T^*_{vert}M)$.

For the general case, let $\alpha \in \Omega^{i,0}(T^*_{vert}M)$ be of
the form $\alpha = \eta \wedge \pi_0^*\beta$, where $\eta \in
\Omega_{\maL}^{i,0}(M)$ and $\beta \in \Omega_c^j(B)$, and $\pi_0 :
T^*_{vert}M \to B$ is the induced projection. Then, using Equation
\eqref{def.*}, we obtain
\begin{equation}
        *_{vert}^2(\alpha) = *_{vert}^2(\eta) \wedge \pi_0^*\beta =
        \alpha.
\end{equation}
Similarly, using the definition of $d_{vert}$, Lemma
\ref{lemma.Bilinear}, and Equation \eqref{def.*}, we obtain
\begin{multline*}
        (-1)^{i+1} d_{vert} \circ *_{vert}(\alpha) = (-1)^{i+1}
        d_{vert} \circ *_{vert}(\eta) \wedge \pi_0^*\beta \\ =
        *_{vert} \circ \delta_{vert}(\eta) \wedge \pi_0^*\beta =
        *_{vert} \circ \delta_{vert}(\alpha).
\end{multline*}
This is enough to complete the proof.
\end{proof}

We are ready now to determine the homogeneous, Laurent-Poisson
homology groups of $A^*(\GR)$. Recall that $p$ denotes the dimension
of the fibers of $M \to B$ and $q$ is the dimension of the manifold
$B$. We set for any fixed $k\in \Z$, 
$$
	K^{j,l} := \maP^{k-j,j}_{l-j},
$$
so that 
$$
	\delta_{vert}:K^{j,l} \to K^{j,l-1} \text{ and } \alpha: K^{j,l} \to K^{j+1,l}.
$$
To compute the homogeneous Laurent $\delta$-homology of
$A^*(\GR)\smallsetminus 0$, we use that the complex splits into
subcomplexes $(\maP^k,\delta)$. Thus we can fix the integer $k\in \Z$  and define a filtration 
of the
above bicomplex $K^{j,l}$ by
$$
	F_{h} := \bigoplus_{l\in \Z, j\leq h} K^{j,l}.
$$

\begin{proposition}\label{delta=delta.vert}\
The spaces $\cohom_{\maL,k+l}^{\delta}(A^*(\GR)\smallsetminus 0\vert
B)_l$ and $\cohom_{\maL,k+l}^{\delta_{vert}}(A^*(\GR)\smallsetminus
0\vert B)_l$ are isomorphic, that is
$$
	\cohom_{\maL,k+l}^{\delta} (A^*(\GR) \smallsetminus 0\vert
	B)_l \simeq \cohom_{\maL,k+l}^{\delta_{vert}} (A^*(\GR)
	\smallsetminus 0 \vert B)_l.
$$
\end{proposition}

\begin{proof}\ Recall that we have 
$$
	\delta = \delta_{vert} + \alpha \text{ and }\delta_{vert}\circ
	\alpha + \alpha \circ \delta_{vert}=0.
$$
We then use for any fixed $k$ the decomposition $\maP^k \cong
\bigoplus_{i+j=k}\maP^{i,j}$ into a finite double complex and the
decreasing filtration $F_{h}$ defined above. This yields a spectral
sequence $(E^r)_{r\geq 1}$ which converges to the $\delta$-homology by
classical homological arguments.  The $E^1$ term of this spectral
sequence is given by
$$
	E^1_{u,v}= \cohom_{\maL}^{2p-v-k+u,u}(T^*_{vert}M
	\smallsetminus 0 \vert B)_{p-k+u}.
$$
But again by a homotopy argument, the homogeneous vertical de Rham cohomology space 
$\cohom_{\maL}^{2p-v-k+u,u}(T^*_{vert}M
\smallsetminus 0 \vert B)_{p-k+u}$ is trivial unless
$u=k-p$. Therefore, we get:
$$
	E^1_{u,v}=0 \text{ if } v\not = -k-p.
$$
Hence for any $r\geq 1$, we see that $d^r=0$ and the spectral sequence
collapses at $E^1$. The proof is thus complete.
\end{proof}

\begin{theorem}\label{theorem.delta}\
The homogeneous, Laurent-Poisson homology groups of $A^*(\GR)$ are
given by
$$
        \cohom_{\maL,k}^{\delta}(A^*(\GR) \smallsetminus 0 \vert B)_l
        \simeq  \cohom_{c,\maL}^{p-l,k-l-p}(T_{vert}^*M
        \smallsetminus 0 \vert B)_0.
$$
\end{theorem}

\begin{proof}\
The vertical symplectic Hodge operator $*_{vert}$ yields
isomorphisms
$$
       	*_{vert} : \bigoplus_{i+j=k}\Omega_{\maL}^{i,j}( T^*_{vert}M
       	\smallsetminus 0 )_l \rightarrow \bigoplus_{i+j=k}
       	\Omega^{2p-i,j}_{\maL}( T_{vert}^*M \smallsetminus 0 )_{l+p-i}
$$
which intertwine the $\delta_{vert}$ and $d_{vert}$ differentials (up
to a sign). Proposition \ref{delta=delta.vert} shows then that
$$
	\cohom_{\maL,k}^{\delta}(A^*(\GR) \smallsetminus 0 \vert B)_l
        \simeq  \bigoplus_{i+j=k}\cohom_{c,\maL}^{2p-i,j}(T_{vert}^*M
        \smallsetminus 0 \vert B)_{l+p-i}.
$$
But for $l+p-i\not =0$, the cohomology spaces
$\cohom^{2p-i,j}(T^*_{vert}M\smallsetminus 0 \vert B)_{l+p-i}$ are
trivial by the homotopy invariance of de Rham cohomology. Hence the
only non trivial term is:
$$
	\cohom_{c,\maL}^{p-l,k-l-p}(T^*_{vert}M\smallsetminus 0 \vert B)_0,
$$
and this completes the proof.
\end{proof}

We can now apply the results of Section~\ref{Sec.filtered} together
with Theorem \ref{theorem.delta}.

\begin{proposition}\label{prop.E2}\
The $\EH^2$-term of the spectral
sequence associated in Lemma \ref{lemma.sp.sq} to the Hochschild homology
of the algebra $\alge_{\maL}(M\vert B;E)$ is given by:
$$
        \EH^2_{k,h} \simeq \cohom_{c,\maL}^{p-k,h-p}
        (S^*_{vert}(M)\times S^1 \vert B),
$$
where $S^*_{vert}(M)$ is the sphere bundle of the vertical cotangent
bundle $T^*_{vert}M$.
\end{proposition}

\begin{proof}\  Denote by $\delta$ the Poisson differential on the
vertically symplectic fibration $T^*_{vert}M \to B$. Lemma
\ref{Lemma.HH.AL}, Proposition \ref{delta=delta.vert} and Theorem
\ref{theorem.delta} give by straightforward computation:
$$
        \EH^2_{k,h} \simeq \cohom_{c,\maL}^{p-k,h-p}(T^*_{vert}M
        \smallsetminus 0 \vert B)_0.
$$
Now a classical argument shows that \cite{MelroseNistor}:
$$
	\cohom_{c,\maL}^{p-k,h-p}(T^*_{vert}M
        \smallsetminus 0 \vert B)_0 \simeq \cohom_{c,\maL}^{p-k,h-p}
        (S^*_{vert}(M)\times S^1 \vert B),
$$
which completes the proof.
\end{proof}

In the following lemma, we shall denote by $\otimes_{tf}$
the completion of the tensor product of two algebras in the
unique natural way that makes the completed tensor product a
topologically filtered algebra. See also Equation \eqref{eq.A3}
where $\otimes_{tf}$ was used before.

\begin{lemma}\label{Lemma.fam.EH}\
Assume that the fibration $\pi:M \to B$ and the bundle of algebras
$\alge_{\maL}(M\vert B;E)$ are trivial; that is, assume that $M =
B\times F$ and $\alge_{\maL}(M\vert B;E) \simeq
\alge_{\maL}(F;E)\otimes_{tf} \CIc(B)$ as topologically filtered
algebras. Then the spectral sequence associated in Lemma
\ref{lemma.sp.sq} to the Hochschild homology of the algebra
$\alge_{\maL}(M\vert B;E)$ collapses at $\EH^2$, converges, and we
have
$$
        \HH_k(\alge_{\maL}(M\vert B;E))\simeq \oplus_{i+j=k}
        \cohom_{c,\maL}^{2p-i} (S^*(F) \times S^1) \otimes
        \Omega_c^j(B).
$$

\end{lemma}

\begin{proof}\
We know from \cite{BenameurNistor} that $\HH_j(\alge_{\maL}(F;E)) \simeq
\cohom_{c}^{2p-j} (\maL(S^*(F))\times S^1)$.  The usual K\"{u}nneth
formula for the Hochschild homology of algebras then gives that
\begin{multline}\label{eq.t.p}
	\HH_k(\CIc(B) \otimes \alge_{\maL}(F;E)) \simeq \oplus
        _{i+j=k} \HH_i(\CI_c(B))\otimes \HH_j(\alge_{\maL}(F;E))\\
        \simeq \oplus_{i+j=k} \Omega_c^i(B) \otimes
        \cohom_{c,\maL}^{2p-j} (S^*(F) \times S^1).
\end{multline}
The last equality is a consequence of the computation of
$\HH_j(\alge_{\maL}(F;E))$ carried out in \cite{BenameurNistor}.
The inclusion $\CIc(B) \otimes \alge_{\maL}(F;E) \to \CIc(B)
\otimes_{tf} \alge_{\maL}(F;E)=\alge_{\maL}(M \vert B; E)$ preserves
the natural filtrations, and hence it induces a morphism of the
corresponding spectral sequences.  By Proposition \ref{prop.E2},
Equation \eqref{eq.t.p}, and naturality, this morphism is an
isomorphism for the $E^2$-terms. We also obtain that all classes in
the $E^2$ term of the spectral sequence associated to $\CIc(B) \otimes
\alge_{\maL}(F;E)$ survive in the $E^\infty$ term (they give rise to
Hochschild homology classes), and hence the spectral sequence
associated to $\CIc(B) \otimes \alge_{\maL}(F;E)$ degenerates at
$E^2$. Consequently, the spectral sequence associated to
$\alge_{\maL}(M \vert B; E)$ also degenerates at $E^2$.

Then Theorem \ref{theorem.conv2} can be used to complete the proof.
\end{proof}

We now extend the above lemma to more general groupoids.

\begin{lemma}\label{lemma.fib.cor}\
The spectral sequence $\EH^r$ associated to the Hochschild homology of
$\alge_{\maL}(M\vert B;E)$ by Lemma \ref{lemma.sp.sq} degenerates at
$\EH^2$ and converges to its Hochschild homology.
\end{lemma}

\begin{proof}\
Denote $\alge = \alge_{\maL}(M\vert B;E)$ in this proof, for
simplicity.  The differential $b$ of the Hochschild complex of $\alge$
is $\CI(B)$-linear, if $f \in \CI(B)$ acts on $a_0 \otimes \ldots
\otimes a_n$ by $f(a_0 \otimes \ldots \otimes a_n) := (fa_0) \otimes
\ldots \otimes a_n$. The filtrations of the Hochschild complex are also
preserved by the multiplication operators with functions $f \in
\CI(B)$. This shows that the spectral sequence associated to the Hochschild homology of
$\alge_{\maL}(M\vert B;E)$ by Lemma \ref{lemma.sp.sq} consists of
$\CI(B)$-modules.

Let $\alge^K$ consist of the complete symbols with support in a fixed
compact set $K \subset M$. Let $K_m \subset M$ be an exhausting
sequence of $M$ by compact subsets (that is, $M = \cup K_m$ and each
$K_m$ is contained in the interior of $K_{m+1}$).  We can assume that
$\alge^{[m]} = \alge^{K_m}$. (This sequence can be reduced to one set if
$M$ is compact.) To prove that $E^2_{k,h}(\alge^{[m]}) = 0$ or that $d^r =
0$ for this spectral sequence, it is enough to check that $f
E^2_{k,h}(\alge^{[m]}) = 0$ (or $fd^r = 0$) for any $f$ with support in
$K_m$. But $f E^2_{k,h}(\alge^{[m]}) = f E^2_{k,h}(\alge)$ if the support
of $f$ is inside $K_m$. This shows that $E^2_{k,h}(\alge^{[m]}) = 0$ if
$E^2_{k,h}(\alge) = 0$.

We now reverse the argument. To prove that $E^2_{k,h}(\alge) = 0$ or
that $d^r=0$ for this spectral sequence, it is enough to check that
$fE^2_{k,h}(\alge) = 0$ (or $fd^r = 0$) for any $f$ in a trivializing
open subset $V$ of $M$ (that is, a subset satisfying the conditions of
our third assumption on the algebra $\alge = \alge_{\maL}(M \vert
B;E)$, Equation \eqref{eq.A3}),  more precisely, satisfying that
\begin{multline}\label{eq.A3'}
        \alge_V := \mathcal
        O(\pi^{-1}(V)) \left ( \Psi^{\infty}(\GR_V;E)/
        \Psi^{-\infty}(\GR_V;E) \right ) \\ \simeq \mathcal
        O(\pi^{-1}(V)) \otimes_{tf}\left ( \Psi^{\infty}(V \times
        \GR_b;E)/ \Psi^{-\infty} (V \times \GR_b;E) \right )
\end{multline}
is a $\CI(B)$-linear isomorphism. Then $fE^2_{k,h}(\alge) = f
E^2_{k,h}(\alge_V)$. By Lemma \ref{Lemma.fam.EH}, $E^r_{k,h}=0$
if $k < -p$. Similarly, $fd^r = 0$ if $r \ge 2$.

This proves that the spectral sequence in Hochschild homology
associated to the algebra $\alge_{\maL}(M\vert B;E)$ by Lemma
\ref{lemma.sp.sq} degenerates at $\EH^2$. It also proves that the
assumptions of Theorem \ref{theorem.conv2} are satisfied, so the
spectral sequence $\EH_{k,h}^r$ converges to Hochschild homology.
\end{proof}

Now we can state and prove the main theorem of this section.
Let $\GR$ be a differentiable groupoid whose space of units is a
manifold with corners $M$ which is the total space of a fibration $\pi
: M \to B$, as before. Denote by $\maF^j$  the local coefficient system defined by
$$
        \maF^{j}(b) := \Hd_{2p-j}(\alge_{\maL} (\pi^{-1}(b))) \simeq
        \cohom_{c,\maL}^{j}(\pi_0^{-1}(b)) \simeq
        \cohom_c^{j}(\maL(\pi_0^{-1}(b)))
$$
with $\pi_0 : S^*_{vert}(M)\times S^1 \to B$ the natural
projection.

\begin{theorem}\label{theorem.fib.cor}\
Assume that $B$ is smooth (without corners) and that $\GR$ satisfies
the assumptions \eqref{eq.A1} and \eqref{eq.A2}. Let
$\alge_{\maL}(M\vert B;E)$ be the algebra of Laurent type complete
symbols on $\GR$ with coefficients in the $\ZZ/2\ZZ$-graded vector
bundle $E$.  Then
\begin{equation}\label{eq.fib.cor}
        \HH_m(\alge_{\maL}(M\vert B;E))\simeq \bigoplus_{k+h=m}
        \Omega_{c}^{h}(B,\maF^{2p-k}).
\end{equation}
\end{theorem}

\begin{proof}\ We apply Lemma \ref{lemma.fib.cor} and deduce that:
$$
       	\HH_m(\alge_{\maL}(M\vert B;E))\simeq \bigoplus_{i+j=m}
       	\EH^2_{i,j}.
$$
The computation of $\EH^2$ was carried out in Proposition \ref{prop.E2} and
the result is:
$$
        \EH^2_{i,j} \simeq \cohom_{c,\maL}^{p-i,j-p}(S^*_{vert}(M)
	\times S^1 \vert B).
$$
On the other hand, Proposition \ref{prop.ind.vert} applied to
$X=S^*_{vert}(M) \times S^1$ gives:
$$
        \cohom_{c,\maL}^{p-i,j-p}(S^*_{vert}(M) \times S^1 \vert B) \simeq
        \Omega_c^{j-p}(B, \maF^{p-i}).
$$
Therefore we get
$$
        \HH_m(\alge_{\maL}(M\vert B;E))\simeq \bigoplus_{i+j=m}
        \Omega_c^{j-p}(B, \maF^{p-i} ).
$$
The conclusion follows by setting $i = k$ and $j = h$.
\end{proof}

\begin{remark}\
The above proposition has to be modified only slightly if $B$ is also
a manifold with corners. For example, when $B$ is compact, the result
remains true if we replace $\Omega^j_c(B)$ with $\mathcal
O(B)\Omega^j(B)$.
\end{remark}

\section{The relative case\label{Sec.Rel}}

Let again $\GR$ be a groupoid with corners satisfying the assumptions
\eqref{eq.A1} and \eqref{eq.A2} with respect to the fibration $\pi:M
\to B$ of the manifold with corners $M$ over the smooth manifold $B$.
Let $X$ be a union of faces of $M$.  We shall denote by $\ideal_X$ the
ideal of smooth functions on $M$ that vanish to infinite order on $X$.

We shall consider in this section the algebra of Laurent complete
symbols on $\GR$ which vanish to infinite order over $X$ and which
represent pseudodifferential operators acting on sections of the
$\ZZ/2\ZZ$-graded vector bundle $E$. This algebra is denoted by
$\alge_{\maL}(M\vert B,X\vert B; E)$. Thus we have:
$$
 	\alge_{\maL}(M\vert B, X\vert B ; E) := \mO \ideal_X
 	\alge(M\vert B; E).
$$
Note that if $X=\emptyset$, then we recover the algebras
$\alge_{\maL}(M\vert B; E)$ studied in the previous sections.  The
proof of Proposition~3 in \cite{BenameurNistor} extends to show that
the algebras $\alge_{\maL}(M\vert B,X\vert B; E)$ are topologically
filtered algebras.

For any fibrations $Y \to M \to B$, we  denote by
$p_Y$ the projection $Y \to M$ and by $\pi_Y$ the composite projection
$Y \to B$. When $Y$ is a manifold with corners, we shall denote by
$p_{\maL}$ the projection $\maL(Y) \to M$ which is the composite map
of $\maL(Y)\to Y$ and $Y\to M$. This last notation is intended to simplify
the  statements of  this section.

Let again $S_{vert}^*(M)$ be the quotient bundle of $T^*_{vert}M
\smallsetminus 0$ by the radial action of $\R^*_+$. In
\cite{BenameurNistor}, the computations of periodic cyclic homology
recover the case of our algebra $\alge_{\maL}(M\vert B, X\vert B ;
E)$.  The result can be stated as follows:

\begin{theorem}\label{periodic2}\cite{BenameurNistor}\
For $q=0,1$, we have:
$$
        \Hp_q(\alge_{\maL}(M\vert B, X \vert B; E)) \simeq
        \cohom_c^{[q]}(\maL(S^*_{vert}M)\times S^1 \smallsetminus
        p_{\maL}^{-1}(X)).
$$
\end{theorem}

\begin{proof}\ Again by a Morita equivalence argument we can forget the
bundle $E$.  We apply Proposition 5 in \cite{BenameurNistor} and
obtain for our groupoid $\GR$:
$$
        \HP_q(\alge_{\maL}(M\vert B, X \vert B)) \simeq \cohom_c^{[q]}
        (\maL(S^*(\GR))\times S^1 \smallsetminus
        {p}_{\maL}^{-1}(X)),
$$
But
$$
        \cohom_c^{[q]}
        (\maL(S^*(\GR))\times S^1 \smallsetminus
        {p}_{\maL}^{-1}(X)) \simeq \cohom_{c,\maL}^{[q]}
        ((A^*(\GR)\smallsetminus 0) \smallsetminus
        p_{A^*(\GR) \smallsetminus 0}^{-1}(X))_0.
$$
{}From Assumption \eqref{eq.A2}, we thus deduce as in Lemma \ref{lemma.is.po}
that:
$$
        \cohom_{c,\maL}^{[q]}
        ((A^*(\GR)\smallsetminus 0) \smallsetminus
        {p}_{A^*\GR\smallsetminus 0}^{-1}(X))_0 \simeq \cohom_{c,\maL}^{[q]}
        ((T_{vert}^*M\smallsetminus 0) \smallsetminus
        p_{T_{vert}^*M\smallsetminus 0}^{-1}(X))_0.
$$
The space $\cohom_{c,\maL}^{[q]}((T_{vert}^*M\smallsetminus 0)
\smallsetminus
p_{T_{vert}^*M\smallsetminus 0}^{-1}(X))_0$ is again isomorphic to the space
$\cohom_c^{[q]}(\maL(S^*_{vert}M)\times S^1 \smallsetminus
p_{\maL}^{-1}(X))$, and this completes the proof.
\end{proof}

Let us state now the corresponding results for Hochschild
homology. Let $\maF_X^l$ be the coefficient system over $B$ given
for any $b\in B$ by the relative cohomology space
$$
     	\maF_X^l(b) := \cohom^l_{\maL}(S^*(\pi^{-1}(b)) \times S^1,
	p_{S^*_{vert}M\times S^1}^{-1}(X)).
$$

\begin{theorem}\
The Hochschild homology spaces of the algebra
$\alge_{\maL}(M\vert B, X \vert B; E)$
are given by:
$$
        \Hd_m(\alge_{\maL}(M\vert B, X \vert B; E)) \simeq
        \oplus_{k+h=m} \Omega_c^h (B, \maF_X^{2p-k}).
$$
\end{theorem}

\begin{proof}\ We can again assume that the  graded vector bundle $E$
is trivial, one-dimensional.  The
proof of this theorem is similar to the proof of Theorem
\ref{theorem.fib.cor}, replacing cohomology by relative cohomology and
using excision. More precisely, the natural filtration of the
topologically filtered algebra $\alge_{\maL}(M\vert B, X \vert B)$ by
the order of the symbols gives rise to a spectral sequence for
Hochschild homology given by Lemma \ref{lemma.sp.sq}. This spectral
sequence was studied in \cite{BenameurNistor} where the $E^2$-term was
identified with the homogeneous Laurent Poisson relative homology of
the Poisson manifold $A^*(\GR)\smallsetminus 0$. Indeed, we have
\cite{BenameurNistor}[Proposition 7]:
$$
        \EH^2_{k,h} \simeq
	\cohom_{\maL,k+h}^{\delta}((A^*(\GR)\smallsetminus 0)
	\smallsetminus p^{-1}_{A^*(\GR)\smallsetminus 0}(X))_k,
$$
 Now using
the immediate extension of Lemma \ref{lemma.is.po} to the relative case, we
obtain:
$$
        \cohom_{\maL,k+h}^{\delta}((A^*(\GR)\smallsetminus 0)
        \smallsetminus p^{-1}_{A^*(\GR)\smallsetminus 0}(X))_k \simeq
        \cohom_{\maL,k+h}^{\delta}((T^*_{vert}M\smallsetminus 0)
        \smallsetminus p^{-1}_{T^*_{vert}M\smallsetminus 0}(X))_k.
$$
The vertical symplectic Hodge operator $*_{vert}$ preserves the
forms vanishing above $X$, therefore
we deduce using the proof of Proposition \ref{prop.E2} that:
$$
        \EH^2_{k,h} \simeq \cohom_{c,\maL}^{p-k,h-p}((S^*_{vert}(M)
	\times S^1) \vert B, p_{\maL}^{-1}(X)\vert B).
$$
The arguments ensuring the degeneracy of the spectral sequence at the
second level in the proof of Lemma \ref{lemma.fib.cor} again obviously
extend to the relative case. In addition, we can apply Theorem
\ref{theorem.conv2} to deduce the convergence of the spectral sequence
to Hochschild homology.  Hence we finally obtain:
\begin{multline*}
        \Hd_m(\alge_{\maL}(M\vert B, X \vert B; E)) \simeq
        \oplus_{k+h=m} \EH^2_{k,h}\\ \simeq \oplus_{k+h=m}
        \cohom_{c,\maL}^{p-k,h-p}((S^*_{vert}(M) \times S^1) \vert B,
        p_{\maL}^{-1}(X)\vert B).
\end{multline*}
To end the proof we simply observe that Proposition
\ref{prop.ind.vert} extends straightforward to the relative case.
\end{proof}

\section{Examples and applications\label{Sec.Examples}}

We begin by providing a construction of a groupoid $\GR$ satisfying
the assumptions \eqref{eq.A1} and \eqref{eq.A2}, for any fibration
$\pi : M \to B$, where $B$ is a smooth manifold (no corners) and $M$
is a manifold, possibly with corners.  For each fibration $\pi$ as
above we shall construct a canonical groupoid $\GR_{M,b}$, the
``$B$-groupoid,'' satisfying the assumptions \eqref{eq.A1} and
\eqref{eq.A2}.

Our examples are obtained by first defining the Lie algebroid of
$\GR$, and then by integrating it. Let $\mathcal V_b(M\vert B)$ be the
space of vertical vector fields on $M$ that are tangent to all faces
of $M$. Then $\mathcal V_b(M\vert B)$ is a Lie algebra with respect to
the Lie bracket of vector fields and is also a projective
$\CI(M)$-module. By the Serre-Swan theorem \cite{Karoubi} there exists
a vector bundle $A_b \to M$ such that $\mathcal V_b(M\vert B) \cong
\Gamma(A)$, naturally. (See also Melrose and Piazza
\cite{MelrosePiazza}.)

The procedure in \cite{NistorINT} then provides us with groupoids
$\GR_{M,b}$ whose Lie algebroid are isomorphic to $A_b$. The minimal
groupoid with this property is obtained as follows. For each face $F$
of $M$, consider the interior of that face, $F_0 := F \smallsetminus
\pa F$. Assume $F$ has codimension $k$. Define then $\GR_F'$ to be the
groupoid with units $F_0$ and with at most one arrow between any two
units. (Thus $\GR_F'$ is the groupoid associated to an equivalence
relation:\ two units are connected by an arrow if, and only if, they
are equivalent.) The equivalence relation that we consider is that $x,
y \in F_0$ are equivalent if, and only if, they belong to the same
connected component of a set of the form $\pi^{-1}(b) \cap F_0$.  Let
\begin{equation}
	\GR_{M,b} = \cup \GR_{F}' \times \RR^k,
\end{equation}
the union being a disjoint union, and with the induced groupoid
structure.  Then it can be checked directly that the charts provided
in \cite{NistorINT} define a smooth structure on $\GR$ such that
$A(\GR_{M,b}) \cong A_b$. This structure must then be unique
\cite{NistorINT}. See also \cite{CrainicFernandes,Monthubert1}.

Let $\mathcal F$ be the locally constant sheaf (or coefficient system)
that associates to $b \in B$ the  complex vector space with basis
the minimal faces of $\pi^{-1}(b)$. (All faces of a manifold with
corners are connected, by definition.)

\begin{theorem}\label{theorem.traces}\ 
Let $\pi : M \to B$ be as above. Then
\begin{equation*}
	\Hd_0(\alge_{\maL}(M \vert B; E)) \cong 
	C_c^{\infty}(B, \mathcal F),
\end{equation*}
the space of compactly supported sections of the sheaf $\mathcal
F$. The space of traces of $\alge_{\maL}(M \vert B; E)$ identifies
with the dual of this space:\
$$
	\Hd^0(\alge_{\maL}(M \vert B; E)) \cong C^{-\infty}(B, \maF)
	=: C_c^{\infty}(B, \mathcal F)' .
$$
\end{theorem}

\begin{proof}\ This follows from Theorem \ref{theorem.fib.cor}.
\end{proof}

In the particular case when $M$ is smooth, our construction simplifies
and we obtain $\GR_{M,b} = M \times_B M$. Then the algebra
$\tPS{\infty}$ consists of differentiable families of
pseudodifferential operators along the fibers of $M \to B$.
Similarly, let us consider $\Psi^\infty(M\vert B;E)$, the algebra of
{\em smooth} families of pseudodifferential operators along the fibers
with coefficients in the $\ZZ/2\ZZ$-graded vector bundle $E$,
introduced in \cite{AtiyahSinger4}.  Let $\alge (M\vert B;E):=
\Psi^{\infty}(M \vert B;E)/\Psi^{-\infty}(M \vert B; E)$ be the
algebra of vertical complete symbols. Then $\alge (M\vert B;E)
\alge_{\maL} (M\vert B;E)$, and we obtain

\begin{corollary}\label{FIB.HH}\
Let $\maF^{j}$ be the locally constant sheaf given by 
the cohomology of the fibers of $S_{vert}^*(M)\times S^1 \to B$.
Then we have:
$$
 	\Hd_m(\alge (M \vert B;E)) \simeq \oplus_{k+h=m}
	\Omega_c^h(B,\maF^{2p-k}).
$$
\end{corollary}

In particular, for $m=0$ and provided that the fibers of $M \to B$ are
connected, this isomorphism becomes
$$
 	\Hd_0(\alge(M\vert B;E))\simeq  \CIc(B).
$$
Therefore, the space of traces is given by
\begin{equation}\label{eq.traces.conn}
	\Hd^0(\alge(M\vert B;E))\simeq \mathcal C^{-\infty}(B) :=
	\CIc(B)',
\end{equation}
the space of distributions on the base manifold $B$.

Let $\omega$ be as before the vertical symplectic form on $T_{vert}^*M
\to B$. The isomorphism of Equation \eqref{eq.traces.conn} can be made
more explicit as follows. Let ${\mathcal R}$ be the radial vector
field on the fibers of $T_{vert}^*M \to B$ and $\alpha=i_{\mathcal
R}(\omega^p/p!)$ the corresponding Liouville form on the fibers of
$S_{vert}^*(M) \to B$.  Let $\STr$ be the graded trace on the endomorphisms
of the fibers of $E$, that is $\STr(A)=\Tr(\gamma \circ A)$, where
$\gamma$ is the involution defining the grading. Also, let $\pi_*$ be
the fiberwise integration on the fibers of $T_{vert}^*M$. Then, for
any $ a=\sum_{j\leq m} a_j$, we set
\begin{equation}
 	\tau_{\mu}(a):= \langle \mu, \pi_*\STr(a_{-p})\alpha \rangle.
\end{equation}
This formula defines a super (or graded) trace on $\alge(M\vert B;E)$
such that $\mu \to \tau_{\mu}$ is the isomorphism described in
Corollary \ref{FIB.HH}, for $m = 0$.

Fix a quantization function 
\begin{equation}
	q : \cup S^s(T^*_{vert}M; \End(E)) \to \cup \Psi^{s}(M \vert
	B;E),
\end{equation}
where $S^s(T^*_{vert}M; \End(E))$ denotes classical vertical symbols of order $s$.
The function $q$ is thus assumed to be continuous and to satisfy
$\sigma_s(q(a)) = a$ if $a$ is a symbol of order $s$.  Let $\rho$ be a
positive symbol on $T^*_{vert}M$, such that $\rho(\xi) = |\xi|$, for
$|\xi| \ge 1$. We consider a family $D(z)$ such that
\begin{equation}
	D(z) = q(\rho^z) B_1(z) + R(z),
\end{equation}
where $B_1$ is a holomorphic function on $\CC$ with values in
$\Psi^{0}(M \vert B;E)$, $B_1(0) = 1$, $R(z)$ is a holomorphic
function on $\CC$ with values $\Psi^{-\infty}(M \vert B;E)$, and
$R(0)=0$. Then $D(0) = 1$.

\begin{proposition}\ Let $\mu$ be a distribution on $B$
and $A \in \Psi^{m}(M \vert B;E)$ and $D(z)$ be as above. Then the
function
$$ 
	z \to F_A(z) :=\langle \mu, \STr_b(A D(z)) \rangle
$$
is well defined for $Re(z) < - m - p$, where $p$ is the dimension of
the fibers of $\pi : M \to B$. The function $F_A$ extends to a
meromorphic function on $\CC$, with at most simple poles at the
integers. For $z = 0$, the residue of this function is up to constant
$\tau_\mu(A)$.
\end{proposition}

\begin{proof}\ This is proved as in the classical case
when $D(z)$ is given by the complex powers of a positive elliptic
operator (see \cite{Seeley,GuilleminDuistermaat}). One can follow the
approach from \cite{NistorBLG}, for example.
\end{proof}

We remark that we used only a weak result from \cite{NistorINT}
(namely Theorem 2).  Let us use this opportunity however to mention
that there is a missing assumption in the gluing theorem of
\cite{NistorINT}, Theorem 3. Here is the corrected version.

\begin{theorem} \label{Theorem.Glueying}\
Let $A$ be a Lie algebroid on a manifold with corners $M$.  Suppose
that $M$ has an $A$-invariant stratification $M = \cup S$ such that,
for each stratum $S$, the restriction $A_S$ is integrabl and let
$\GR_S$ be $d$-simply connected differential groupoids such that
$A(\GR_S) \simeq A_S$. Then $A$ is integrable if, and only if, the
exponential map $\Exp : A \to \GR = \cup \GR_S$ is injective on an
open neighborhood of the zero section of $A$ for some (equivalently,
for any) connection on $A$. Moreover, if these conditions are
satisfied, then the disjoint union $\GR = \cup \GR_S$ is naturally a
differentiable groupoid such that $A(\GR) \simeq A$.
\end{theorem}

The above condition on the map $\Exp$ is seen to be necessary from the
work of Crainic and Fernandes \cite{CrainicFernandes}, and also from
some earlier results of Weinstein. The map $\Exp$ introduced in
\cite{NistorINT} seems to be essential for both the results of that
paper and for the results of \cite{CrainicFernandes}. The second named
author would like to thank M. Crainic for pointing out a possible
problem with the original statement of the above theorem.

Let us mention for completeness the result for the computation of the
cyclic homology of the algebras $\alge_{\maL}(M\vert B;E)$.  Note
first that
\begin{equation*}
	\Hc_j(\alge_{\maL}(M\vert B;E)) \simeq
	\Hp_j(\alge_{\maL}(M\vert B;E)) \quad \text{for } j \ge 2p +
	2,
\end{equation*}
from the SBI-exact sequence, because the Hochschild homology of
$\alge_{\maL}(M\vert B;E)$ vanishes above this rank.

\begin{proposition}\label{prop.cyclic.A} 
The spectral sequence $\EC^r_{k,h}$ associated to the cyclic homology
complex of $\alge_{\maL}(M\vert B;E)$ has $E^1$-term given by $\EC^1_{k,h} =
\Omega_{rc,\maL}^{k+h}(S_{vert}^*M)$, if $k \not = 0$, and by 
$$
	\EC^1_{0,h} =
	\Omega^h(S_{vert}^*M \times S^1)/ d\Omega^{h-1}(S_{vert}^*M \times S^1)
	\oplus \oplus_{j > 0}\cohom^{h-2j}(S^*_{vert}M \times S^1).
$$
The $d_1$ differential is induced by the Poisson diff\'erential and
the spectral sequence converges to the cyclic homology of the algebra
$\alge_{\maL}(M\vert B;E)$.
\end{proposition}

%MEMO: 
%AN INTERESTING PROBLEM HERE FOR FUTURE: ?HOMOGENEOUS HOMOLOGY OF
%$(\Omega/d\Omega,\delta)$ NEEDS TO BE LINKED WITH SOME CONTACT
%HOMOLOGY. VICTOR, I'M WONDERING IF SOMEONE LIKE WEINSTEIN CAN ANSWER
%THIS KIND OF QUESTION?


\begin{thebibliography}{10}
\bibitem{AtiyahSinger4} M. {F}. Atiyah and {I}. {M} Singer, \emph{The
index of elliptic operators IV}, {A}nn. of Math. \textbf{93} (1971),
119--138.

\bibitem{AtiyahMacDonald}
M. F. Atiyah and MacDonald, \emph{Introduction to commutative algebra},
Addison-Wesley, Reading, Mass.-London, 1969.

\bibitem{AtiyahBottPatodi} M.{F}. Atiyah, R. Bott, and V. Patodi,
\emph{On the heat equation and the index theorem},
Inv. Math. \textbf{19}
(1973), 279--330.

\bibitem{BenameurNistor} M.-{T}. Benameur and {V}. {N}istor,
\emph{Homology of complete symbols and Noncommutative geometry}, to
appear in "Collected papers on Quantization of Singular Symplectic
Quotients", Progress in Math. series of Birkh\"{a}user, to appear.

\bibitem{BenameurNistorFol} M.-{T}. Benameur and {V}. {N}istor,
\emph{Residues and an index theorem for foliations}, work in progress.

\bibitem{BerlineGetzlerVergne} N.~Berline, E. Getzler and M. Vergne,
\emph{Heat kernels and Dirac operators}, A series of Comprehensive
studies in Mathematics, Springer-Verlag 1991.

\bibitem{Bismut} J.-L. Bismut, \emph{The Index Theorem for families of
Dirac operators: two heat equation proofs}, Invent. Math. \textbf{83}
(1986), 91--151.

\bibitem{BismutCheeger} J. M. Bismut and J. Cheeger, {\em Families
index for manifolds with boundary, superconnections, and cones.
I. Families of manifolds with boundary and Dirac operators.},
J. Funct.  Anal. 89 (1990), 313--363.

\bibitem{Brylinski} J.-{L}. Brylinski, \emph{A differential complex
for {P}oisson manifolds}, J.{D}iff. {G}eom.  \textbf{28} (1988),
93--114.

%\bibitem{BenameurNistor} M. Benameur and V. Nistor, \emph{Residues and an
%index theorem for
%foliations},  work in progress.

\bibitem{BrylinskiGetzler} J.-{L}. Brylinski and {E}. {G}etzler,
\emph{The homology of {A}lgebras of {P}seudo-differential {S}ymbols
and the {N}oncommutative {R}esidue}, $K$-{T}heory \textbf{1} (1987),
385--403.

\bibitem{ConnesNCG} A.~Connes, \emph{Noncommutative differential
geometry}, Publ. Math. IHES \textbf{62} (1985), 41--144.

\bibitem{ConnesBOOK} A.~Connes.  \newblock {\em Noncommutative
Geometry}, \newblock Academic Press, New York - London, 1994.

\bibitem{ConnesMoscovici} A.~Connes and H. Moscovici, \emph{The local
index formula in noncommutative geometry}, Geom. Funct. An. Vol
\textbf{5}, No 2, (1995).

\bibitem{CrainicFernandes} M. Crainic and R. L. Fernandes, Integrability of
Lie brackets, math.DG/0105033.

\bibitem{Getzler} E.~Getzler, \emph{Pseudodifferential operators on
supermanifolds and the Atiyah-Singer index theorem},
Comm. Math. Phys. \textbf{92}, (1983) , 163-178.

\bibitem{Guillemin} V.~Guillemin, \emph{A new proof of Weyl's formula
on the asymptotic distribution of eigenvalues},
Adv. in Math. \textbf{55}, (1985) , 131-160.

\bibitem{GuilleminDuistermaat} J.J Duistermaat and V.~Guillemin,
\emph{ The spectrum of positive elliptic operators and periodic
bicharacteristics}, Invent. Math. \textbf{29}, (1975), no. 1, 39--79.

\bibitem{Kassel} C. Kassel, \emph{A K\"{u}nneth formula for the cyclic
cohomology of ${Z}/2$-graded algebras}, Math. Ann. 275 (1986), no. 4,
683--699.

\bibitem{KaroubiCyclic} M. Karoubi, \emph{Formule de Kunneth en
homologie cyclique.}, C. R. Acad. Sci. Paris S\'er. I Math. 303 (1986)
no. 13, 595--598.

\bibitem{Karoubi} M. Karoubi, \emph{Homologie cyclique et
$K$-th\'eorie}, Asterisqu\'e, \textbf{149} (198), 147 pp.

\bibitem{LMN} R. Lauter, B. Monthubert, and V. Nistor,
\emph{Pseudodifferential operators on continuous family groupoids},
Doc. Math. (electronic) \textbf{5} (2000) 625--655.

\bibitem{LauterMoroianu} R. Lauter and S. Moroianu, \emph{Homology of
pseudo-differential operators on manifolds with fibered boundaries},
Mainz University preprint 2000.

\bibitem{LauterNistor} R. Lauter and V. Nistor, \emph{Analysis of
geometric operators on open manifolds: a groupoid approach}, to appear in
"Collected papers on Quantization of Singular Symplectic  Quotients",
Progress in Math. series of Birkh\"{a}user, to appear.

\bibitem{LauterSeiler} R.~Lauter and J.~Seiler,
\emph{Pseudodifferential analysis on manifolds with boundary -- a
comparison of b-calculus and cone algebra}, to appear in: J. B. Gil et
al.\ (eds), {\em Approaches to Singular Analysis}, Operator Theory:
Advances and Applications, Birkh{\"a}user, Basel.

\bibitem{LP} M. Lesch and M. J. Pflaum, \emph{Traces on algebras of
parameter dependent pseudodifferential operators and the
eta-invariant}, to appear.

\bibitem{Loday} J-L Loday, \emph{Cyclic homology}, A Series of
comprehensive studies in Mathematics 301, (1992).

\bibitem{Moroianu} S. Moroianu, \emph{Residue functionals on the
algebra of adiabatic pseudo-differential operators}, MIT thesis, 1999.

\bibitem{Loday-Quillen1} J.-L. Loday and D.~Quillen, \emph{Cyclic
homology and the {L}ie homology of matrices},
Comment. Math. Helv. \textbf{59} (1984), 565--591.

\bibitem{Mackenzie} 
K.~Mackenzie, \emph{Lie groupoids and {L}ie algebroids in differential 
geometry}, Lecture Notes Series 124, London Mathematical Society, 1987. 
 
\bibitem{MacLane} S.~Mac~Lane, \emph{Homology}, Springer-Verlag,
Berlin-Heidelberg-New York, 1995.

\bibitem{MacLaneMoerdijk} S.~Mac~Lane and I. Moerdijk, \emph{Sheaves in
geometry and logic. A first introduction to topos theory}, Universitext,
Springer-Verlag, Berlin-Heidelberg-New York, 1994.

\bibitem{MelroseNistor} R. Melrose and {V}. {N}istor, \emph{Homology
of pseudodifferential operators I. Manifolds with boundary}, accepted
for publication in Amer. J. Math.

\bibitem{MelrosePiazza} R. Melrose and {P}. Piazza, \emph{Analytic
K-theory on manifolds with corners}, Adv. Math. 92 (1992), no. 1,
1--26.

\bibitem{Monthubert1} B.~Monthubert, \emph{Pseudodifferential calculus
on manifolds with corners and groupoids},
Proc. Amer. Math. Soc. \textbf{12} (1999), 2871--2881.

\bibitem{Monnier} P. Monnier, \emph{Poisson cohomology in dimension
two}, preprint math.DG/0005261.

\bibitem{NazaShatSter} V. E. Nazaikinskii, V. E. Shatalov, and
B. Yu. Sternin, \emph{Methods of Noncommutative Analysis},
de Gruyter Studies in Mathematics 22, 1996.

\bibitem{Schrohe} E. Schrohe, \emph{Non-commutative residue and
manifolds with conical singularities}, J. Funct. Anal. \textbf{150}
(1997), 146--174.

\bibitem{SchroheN} R. Nest and E. Schrohe, \emph{Hochschild
homology of Boutet de Montvel's algebra}, work in progress.

\bibitem{NistorINT} V. Nistor, \emph{Groupoids and integration
of Lie algebroids}, J. Math. Soc. Japan \textbf{52} (2000), 847--868.

\bibitem{NistorBLG} V. Nistor, \emph{An index theorem for families
invariant with respect to a bundle of Lie groups}, Preprint (2000).

\bibitem{NWX} V.~Nistor, A. Weinstein and Ping Xu,
\emph{Pseudodifferential operators on differential groupoids}, Pacific
J. Math. \textbf{189} (1999), 117-152.

\bibitem{SchulzeSterninShatalov} B.-W. Schulze, V. E. Shatalov, and
B. Sternin, Differential equations on singular 
manifolds. \emph{Semiclassical theory and operator algebras.} Mathematical 
Topics, 15. Wiley-VCH Verlag Berlin Gmbh,
Berlin, 1998, 376 pp.

\bibitem{Shubin} M.~A. Shubin, \emph{Pseudodifferential operators and
spectral theory}, Springer Verlag, Berlin-Heidelberg-New York, 1987.

\bibitem{Seeley} R. T. Seeley, \emph{Complex powers of an elliptic operator},
Proc. Symp. Pure Math. 10, (1967), 288-307.

\bibitem{Taylor} M. Taylor,  \emph{Pseudo differential operators},
Lecture Notes in Mathematics, Vol. 416. Springer-Verlag, Berlin-New York, 1974.

\bibitem{Tsygan} B.~L. Tsygan, \emph{Homology of matrix {L}ie algebras
over rings and {H}ochschild homology.} Uspekhi Math. Nauk. \textbf{38}
(1983), 217--218.

\bibitem{Tondeur} P. Tondeur, \emph{Geometry of foliations.}
Monographs in Math. Birkhauser, 1997.

\bibitem{VaismanBook} I.~Vaisman, \emph{Lecture on the Geometry of
Poisson Manifolds}, Progress in Mathematics, Birkhauser, 1994.

\bibitem{Vaisman2} I.~Vaisman, \emph{New examples of twisted
cohomologies}, Boll. U.M.I.  \textbf{7} (1993), 355--368.

\bibitem{Wodzicki} M.~Wodzicki, \emph{Excision in cyclic homology and
in rational algebraic {K}-theory}.  Annals of Mathematics,
129:591--640, 1989.

\bibitem{WodzickiU} M.~Wodzicki, an unpublished hand written preprint
dated 1989.

\bibitem{Weinstein} A. Weinstein, \emph{The local structure of Poisson
manifolds}, Journal of Diff. Geom. 18 (1983), no. 3, 523--557.

\bibitem{WeinsteinXu} A. Weinstein and Ping Xu, Hochschild cohomology
and characteristic classes for star-products. \emph{Geometry of
differntial equations}, 177--194, Amer. Math. Soc. Transl. Ser 2, 186,
AMS, Providence, RI, 1998.

\bibitem{PingXu} Ping Xu, \emph{Poisson cohomology of regular Poisson
manifolds}, Ann. Inst. Fourier \textbf{42} (1992), 967--988.
\end{thebibliography}
\end{document}